\documentclass[12pt,a4paper,reqno]{amsart}
\usepackage[english]{babel}
\usepackage{amssymb,latexsym,amsfonts,amsthm,upref,amsmath}
\usepackage[margin=1in]{geometry}
\usepackage{hyperref}
\hypersetup{colorlinks,citecolor=blue,filecolor=black,linkcolor=blue,urlcolor=blue}
\usepackage{comment} 
\usepackage{enumerate} 
\usepackage{tikz}
\usepackage{float}
\usepackage{cleveref}

\newtheorem*{thm*}{Theorem}
\newtheorem*{lem*}{Lemma}
\newtheoremstyle{prim}{}{}{\normalfont}{}{\bfseries}{}{ }{}
\theoremstyle{prim}
\newtheorem{ex}{Example}
\newtheoremstyle{stil}{}{}{\slshape}{}{\bfseries}{}{ }{}
\theoremstyle{stil}
\newtheorem{thm}{Theorem}[section]
\newtheoremstyle{defi}{}{}{}{}{\bfseries}{}{ }{}
\theoremstyle{defi}
\newtheorem{defn}[thm]{Definition}
\theoremstyle{defi}
\newtheorem{rem}[thm]{Remark}
\theoremstyle{stil}
\newtheorem{pro}[thm]{Proposition}
\theoremstyle{stil}
\newtheorem{lem}[thm]{Lemma}
\theoremstyle{stil}
\newtheorem{kor}[thm]{Corollary}
\newenvironment{prf}{\noindent \textit{Proof.}}{\null\hfill$\qed$\hskip
2mm\vskip 2mm}

\newcommand{\ch}{\mathop{\mathrm{ch}}}

\newcommand{\q}{{\protect\underline{\mathsf{q}}}}
\newcommand{\ndo}{\mathop{\mathrm{End}}}
\newcommand{\om}{\mathop{\mathrm{Hom}}}

\newcommand{\sym}{\mathop{\mathrm{Sym}}}
\newcommand{\gauss}[2]{\genfrac{[}{]}{0pt}{}{#1}{#2}}

\numberwithin{equation}{section}

\begin{document}

\title[Vertex operators and principal subspaces of level one for $U_q (\widehat{\mathfrak{sl}}_2)$]{Vertex operators and principal subspaces of level one for $U_q (\widehat{\mathfrak{sl}}_2)$}

\author{Slaven Ko\v{z}i\'{c}} 

\address{School of Mathematics and Statistics F07, University of Sydney, NSW 2006, Australia}

\email{kslaven@maths.usyd.edu.au}



\begin{abstract}
We consider two different methods of associating vertex algebraic structures with the level
$1$ principal subspaces for $U_q (\widehat{\mathfrak{sl}}_2)$. 
In the first approach, we introduce certain commutative operators and study the corresponding vertex algebra and its module. We find combinatorial bases for these objects and show that they coincide with the principal subspace bases found by 
B. L. Feigin and A. V. Stoyanovsky. In the second approach, we introduce the, so-called nonlocal $\underline{\mathsf{q}}$-vertex algebras, investigate their properties
and construct the nonlocal $\underline{\mathsf{q}}$-vertex algebra and its module, generated by Frenkel-Jing operator and Koyama's operator respectively. 
By finding the combinatorial bases of their suitably defined subspaces, we establish a connection with the sum sides of the Rogers-Ramanujan identities.
Finally, we discuss further applications  to quantum quasi-particle relations.
\end{abstract}

 \maketitle

\section*{Introduction} 
In \cite{FS}, B. L. Feigin and A. V. Stoyanovsky associated with every standard  $\widehat{\mathfrak{sl}}_2$-module $L$ the principal subspace -- a  subspace of 
$L$ obtained by the action of the universal enveloping algebra $U(\widehat{\mathfrak{n}}_+)$ on the highest weight vector of $L$, where  $\widehat{\mathfrak{n}}_+$ is
a subalgebra of $\widehat{\mathfrak{sl}}_2$ 
corresponding to a certain triangular decomposition. The authors constructed  monomial bases for the principal subspaces and computed the character formulas, which are related to the sum sides of the Rogers-Ramanujan-type identities. This led to an intensive research of  the principal subspaces  and related structures.

Generalizing their work, G. Georgiev in \cite{G} constructed similar bases for the principal subspaces corresponding to certain standard $\widehat{\mathfrak{sl}}_{n+1}$-modules. More recently, the problem of finding such combinatorial bases in types $(BC)^{(1)}$ was solved by M. Butorac (\cite{Bu1}--\cite{Bu2}) for certain highest weights. In \cite{CLM1}--\cite{CLM2}, S. Capparelli, J. Lepowsky and A. Milas applied the theory of vertex algebras and intertwining operators to recover Feigin-Stoyanovsky's character formulas for principal subspaces associated to $\widehat{\mathfrak{sl}}_2$ and to show that the corresponding graded dimensions satisfy Rogers-Selberg recursions. Next, in the papers of C. Calinescu, Lepowsky, Milas and C. Sadowski (\cite{x1}--\cite{x4},\cite{S1}), natural presentations of certain principal subspaces associated with the affine Lie algebras  of types $(ADE)^{(1)}$ and $A_{2}^{(2)}$ were found. 

Another variant of principal subspaces, the so-called Feigin-Stoyanovsky's type subspaces, was studied by M. Primc. 
He constructed monomial bases for the subspaces corresponding to the affine Lie algebras of type $A_{n}^{(1)}$ (\cite{P1}) and, for the highest weight $\Lambda_0$, to the affine Lie algebras of type $(ABCD)^{(1)}$  (\cite{P2}).
Later on, combinatorial bases and character recurrence relations for Feigin-Stoyanovsky's type subspaces were studied by I. Baranovi\'{c} (\cite{Bar}), M. Jerkovi\'{c} and Primc (\cite{J1}--\cite{J2},\cite{JP}) and G. Trup\v{c}evi\'{c} (\cite{T1}--\cite{T3}). 
For more details on the principal subspaces the reader may  consult, for example, the papers \cite{AKS}, \cite{FFJMM}, \cite{Kawa}, \cite{MiP}, \cite{S2} and references therein. 

Using  Drinfeld realization of quantum affine algebras (\cite{D}) and Frenkel-Jing realization of the integrable highest weight modules (\cite{FJ}),  we studied in \cite{Kozic} the notion of a principal subspace associated with the integrable highest weight module for the quantum affine algebra $U_q (\widehat{\mathfrak{sl}}_{n+1})$, which was originally introduced in \cite{DF}. We found, for certain highest weights, combinatorial bases of these subspaces, which were, as in \cite{G}, expressed in terms of monomials of the so-called quasi-particles acting on the highest weight vector. Two main ingredients in this construction were finding the relations among quasi-particles and proving the linear independence by using the explicit formulas for certain vertex operators $\mathcal{Y}(z)$ found by Y. Koyama in \cite{Koyama}.
However, unlike the classical case, these relations did not seem to have any vertex algebraic interpretation.
Although the considered operators, i.e. quasi-particles, are not local, they do have some  properties, such as $q$-integrability (see \cite{DM}), which might suggest the existence of an underlying vertex algebra-like theory.

By now, there were several approaches to the development of the quantum vertex algebra theory, which led to the construction of some new and important structures:
$H_D$-quantum vertex algebras of I. I. Anguelova and M. J. Bergvelt (\cite{AB}), 
field algebras of B. Bakalov and V. G. Kac (\cite{Bak}), (different notions of) quantum vertex algebras
of R. E. Borcherds (\cite{Bor2}) and H.-S. Li (\cite{LiNonlocal}--\cite{Li3}), quantum VOAs of P. Etingof and D. Kazhdan (\cite{EK}) and deformed chiral algebras of E. Frenkel and N. Reshetikhin (\cite{FR}). Even though we found in \cite{Kozic} an application of Li's theory  on the quasi-particle construction, we were not able to further develop
this correspondence. In this paper we study two different methods of associating vertex algebraic  structures with the level 1 principal subspaces of $U_q (\widehat{\mathfrak{sl}}_2)$.

Denote by $L=L_0 \oplus L_1$ a direct sum of (nonequivalent) level 1 integrable highest weight $U_q (\widehat{\mathfrak{sl}}_2)$-modules $L_0$, $L_1$ with integral dominant highest weights $\Lambda_0$, $\Lambda_1$ respectively. In  \Cref{komutativnisection}, motivated by Ding-Feigin construction of commutative operators in \cite{DF}, we introduce certain new operators on $L$, 
$\widehat{x}(z)\in\om(L,L((z)))$ and $\widehat{\mathcal{Y}}(z)\in \om(L,L((z^{1/2})))$. These operators were designed to satisfy the following equalities:
\begin{align}
&[\widehat{x}(z_1),\widehat{x}(z_2)]=[\widehat{x}(z_1),\widehat{\mathcal{Y}}(z_2)]=0,\label{int1}\\ &\widehat{x}(z)\widehat{x}(z)=\widehat{x}(z)\widehat{x}'(z)=\widehat{x}(z)\widehat{\mathcal{Y}}(z)=0.\label{int2}
\end{align}
Furthermore, the classical limit of $\widehat{x}(z)$ coincides with the classical limit of the corresponding Frenkel-Jing operator $x^{+}(z)$ (see \cite{FJ}).
Commutativity condition \eqref{int1} will allow us to easily apply the vertex algebra theory, while the second condition, \eqref{int2} will be crucial in finding (independent) relations among the operators. Even though  Ding-Feigin operators $\bar{x}^{+}(z)$ have similar properties, their $q$-integrability relations (cf. \cite{DF}),
$\bar{x}^{+}(z)\bar{x}^{+}(zq^{\pm 2})=0$, seem to be much harder to handle by vertex algebraic methods, than the (classic) integrability relation \eqref{int2} (cf. \cite{LP}). 

First, we construct bases $\widehat{\mathfrak{B}}_{W(\Lambda_i)}$, $i=1,2$, of the principal subspaces $W(\Lambda_i)$ associated with $L_i$. The basis elements are expressed as monomials
of the operator's $\widehat{x}(z)=\sum_{r\in\mathbb{Z}}\widehat{x}(r)z^{-r-1}$ coefficients acting on the highest weight vector $v_{\Lambda_i}$, whose indices satisfy certain difference and initial conditions. These conditions coincide with the conditions found in \cite{G}, while the proof relies upon the results in \cite{Kozic}.
Next, we consider vertex algebra $W_0\subset\om(L_0,L_0 ((z)))$, generated by (a restriction of) $\widehat{x}(z)$, and $W_0$-module $W_1 \subset\om(L_0,L_1 ((z)))$, generated by (a restriction of) $\widehat{\mathcal{Y}}(z)$. By using relations \eqref{int1}, \eqref{int2} we find bases $\widehat{\mathfrak{B}}_{W_i}$, $i=0,1$, of $W_i$:
\begin{align*}
\widehat{\mathfrak{B}}_{W_i}=\bigg\{ \widehat{x}(z)_{l_m}\ldots\widehat{x}(z)_{l_2}\widehat{x}(z)_{l_1}a_i (z) \,\,\bigg| \,\,\,\bigg.
&l_1\leq -1-i\text{ and }\bigg.
 l_{r+1}\leq l_r -2\\
\bigg. &\text{for all }l_r\in\mathbb{Z},\,r=1,2,\ldots,m-1,\,
m\in\mathbb{Z}_{\geq 0}\bigg\},\nonumber
\end{align*}
where $a_0(z)=1$ and $a_1 (z) =\widehat{\mathcal{Y}}(z)$. The form of the both pairs of bases, $\widehat{\mathfrak{B}}_{W(\Lambda_i)}$ and $\widehat{\mathfrak{B}}_{W_i}$, allows us to easily see that the characters of $W(\Lambda_i)$ and $W_i$ are equal to the sum sides of the famous Rogers-Ramanujan identities (cf. \cite{A1}),
\begin{align*}
\prod_{r \geq 0}\frac{1}{(1-q^{5r+1+i})(1-q^{5r+4-i})}&= \sum_{r\geq 0}\frac{q^{r^2+ir}}{(1-q)(1-q^{2})\cdots (1-q^r)}.
\end{align*}

In the second approach, we work with   Frenkel-Jing operator $x(z)=x^{+}(z)$ acting on $L$.
Since $x(z)$ is not local, our first step is  developing an appropriate vertex algebraic setting, which would allow us
to apply similar techniques as in  \Cref{komutativnisection}.
In \Cref{teorijskisection}, we study the operators $a_1 (z),..., a_n(z)\in\mathcal{E}(L)=\om(L,L((z)))$ satisfying the following, so-called quasi-commutativity property:
$$a_1(z_1)a_2(z_2)\cdots a_n(z_n) \in \om(L,L((z_1,z_2,...,z_n))),$$
which is a special case of quasi-compatibility (cf. \cite{LiNonlocal}).
Motivated by \cite{Bak} and  \cite{LiNonlocal}, we introduce the notion of nonlocal $\q$-vertex algebra,  obtained from the notion of vertex algebra (cf. \cite{LiLep}) by replacing
Jacobi identity with ``$\q$-associativity", 
\begin{equation}\label{intro_as}
Y(a(z),z_0 +z_2)Y(b(z),z_2) c(z)=Y(Y(a(z),z_0)b(z),z_2)c(z),
\end{equation}
where variables $z_2 ,z_0,z$ satisfy the following noncommutative constraints
$$z_2 z_0=\q z_0 z_2,\qquad z_2 z=\q zz_2,\quad z_0 z=\q z_0 z$$
for some transcendental element $\q$ over $\mathbb{C}$.
In general, the above definition may be given in a more abstract form, written in terms of {\em states}, instead of in terms of {\em fields}.
However, in this paper, \eqref{intro_as} was easier to handle and we did not need a more general definition. Even though the operators studied in this paper satisfy ``$\q$-associativity'',  we expect that, in a more general setting, a weaker form of this axiom should be considered. 
Finally, we introduce  $r$th products, $r\in\mathbb{Z}$, among quasi-commutative operators in $\mathcal{E}(L)$ and prove the main result of this section:

\makeatletter
\def\thmhead@plain#1#2#3{%
  \thmname{#1}\thmnumber{\@ifnotempty{#1}{ }\@upn{#2}}%
  \thmnote{ {\the\thm@notefont#3}}}
\let\thmhead\thmhead@plain
\makeatother  

\begin{thm*}[\textbf{\ref{main}}]
Let $\mathcal{S}$ be a quasi-commutative subset of $\mathcal{E}(L)$. There exists a unique
smallest nonlocal $\q$-vertex algebra
$V\subset \mathcal{E}(L)$ such that $\mathcal{S}\subseteq V$.
\end{thm*}

In \Cref{FJsection},  using Theorem \ref{main}  we construct nonlocal $\q$-vertex algebra $\left<x(z)\right>$ (for $\q=q^2$), generated by $x(z)$, and a $\left<x(z)\right>$-module $\left<\mathcal{Y}(z)\right>$, generated by $\mathcal{Y}(z)$. The underlying vector spaces, obtained in this way, 
are ``much bigger'' than their analogues from   \Cref{komutativnisection}, so we consider only their subspaces 
$W_{i,\q}$, $i=0,1$, spanned by the operators
\begin{equation}\label{intro_op}
 x(z)_{l_m}\ldots x(z)_{l_1}a_{i}(z),\quad l_j \leq -1,\, j=1,2...,m,\, m\in\mathbb{Z}_{\geq 0},
\end{equation}
where $a_0 (z) =1$, $a_{1} (z) = \mathcal{Y}(z)$. More precisely, the operators
\begin{equation*}
 \widehat{x} (z)_{l_m}\ldots \widehat{x} (z)_{l_1}\widehat{a}_{i}(z),\quad l_j \leq -1,\, j=1,2...,m,\, m\in\mathbb{Z}_{\geq 0},
\end{equation*}
span the whole $W_i$, while the operators \eqref{intro_op} span only the subspace $W_{i,\q}$ of  $\left<x(z)\right>$ (when $i=0$) or 
$\left<\mathcal{Y}(z)\right>$ (when $i=1$). This was caused by the application of the two crucially different approaches. While in 
\Cref{komutativnisection} we adjusted the operators $\widehat{x}(z)$ and $\widehat{\mathcal{Y}}(z)$, so that they can be efficiently handled by vertex algebraic methods,
in Sections \ref{teorijskisection} and \ref{FJsection}, we adjusted the (nonlocal) vertex algebra theory, so that it can handle the original operators $x(z)$ and 
$\mathcal{Y}(z)$. In the end, $r$th products among quasi-commutative operators, which were designed to match the $q$-integrability satisfied by Frenkel-Jing operators,
$x(z)x(z)=x(z)x(zq^2)=0$, gave rise to a ``bigger'' structure, than the usual vertex operator products among local operators satisfying \eqref{int1}.
Since the definition of the subspaces $W_{i,\q}$, given in \eqref{intro_op}, corresponds with the original definition of the principal subspaces, we proceed to study these subspaces.

The main result in \Cref{FJsection} is the construction of the monomial bases $\mathfrak{B}_{W_{i,\q}}$ for the spaces $W_{i,\q}$,
 \begin{align*}
\mathfrak{B}_{W_{i,\q}}=\bigg\{ x(z)_{l_m}\ldots x(z)_{l_1}a_i (z) \,\,\bigg| \,\,\,\bigg.
&l_1\leq -1-i\text{ and }\bigg.
 l_{r}\leq -3\\
\bigg. &\text{for all }l_r\in\mathbb{Z},\,r=2,3,\ldots,m,\,
m\in\mathbb{Z}_{\geq 0}\bigg\}.\nonumber
\end{align*} Using the nonlocal $\q$-vertex algebra setting we  find a construction which is to a great extent analogous to the construction of bases  $\widehat{\mathfrak{B}}_{W_i}$. 
These bases are not given in terms of the same different and initial conditions as the bases $\widehat{\mathfrak{B}}_{W_i}$. 
However, by constructing a bijection between the sets of diagrams corresponding to the elements of $\mathfrak{B}_{W_{i,\q}}$ and $\widehat{\mathfrak{B}}_{W_i}\equiv\widehat{\mathfrak{B}}_{W(\Lambda_i)}$,
we show that, for a suitably defined character $ch_{\q}$, we get
\begin{thm*}[\textbf{\ref{rr}}]
For $i=0,1$ we have
 $$\textstyle\ch_{\q}\displaystyle W_{i,\q}=\sum_{r\geq 0}\frac{q^{r^2+ir}}{(1-q)(1-q^{2})\cdots (1-q^r)}.$$
\end{thm*}

In \Cref{zadnjisection}, we recall the notions of quasi-particles from \cite{G} and (quantum) quasi-particles from \cite{Kozic}, which were, roughly speaking, the
main building blocks of the combinatorial bases found in these papers. One of the most important ingredients in such constructions
is finding an appropriate set of relations among quasi-particles. Motivated by the vertex algebraic interpretation of the quasi-particle relations found by Georgiev
for $\widehat{\mathfrak{sl}}_{n+1}$,
as an application of nonlocal $\q$-vertex algebras we find a similar interpretation of the quasi-particle relations for $U_q (\widehat{\mathfrak{sl}}_{n+1})$.

\allowdisplaybreaks

\section{Preliminaries}

\subsection{Quantum calculus}
We briefly recall some elementary notions of quantum calculus.
For more details the reader may consult \cite{qKac}.
Fix an indeterminate $\q$. First, we define some elements of
the field $\mathbb{C}(\q)$.
For any two integers $m$ and
$l$, $l\geq 0$, define {\em$\q$-integers}, $$[m]_{\q}=\frac{\q^m -1}{\q-1}=1+\q+...+\q^{m-1}, $$
 {\em $\q$-factorials},
 $$[0]_{\q} !=1,\quad [l+1]_{\q}!=[l+1]_{\q}[l]_{\q}\cdots[1]_{\q},$$
 and {\em$\q$-binomial coefficients},
 $$\gauss{m}{l}_{\q}=\frac{[m]_{\q}[m-1]_{\q}\cdots [m-l+1]_{\q}}{[l]_{\q}!}.$$
 
Denote by $z_0$ and $z$  two
noncommutative variables satisfying 
\begin{equation}\label{wz}
z_0 z=\q zz_0.
\end{equation}

\begin{rem}
In the rest of this paper we shall assume that all  formal variables are commutative, unless stated otherwise (as above).
\end{rem}

We have the following $\q$-analogue of the binomial theorem:
 \begin{pro}\label{expand}
 For every integer $m$ and variables $z_0 ,z$ satisfying \eqref{wz} we have
 \begin{equation}\label{qbinomial}
(z + z_0 )^m=\sum_{l\geq 0} \gauss{m}{l}_{\q} z^{m-l}z_{0}^{l}.
\end{equation}
 \end{pro}
 
 Let $V$ be a vector space over the field $\mathbb{C}(\q)$ and let
 $a(z)\in V[[z^{\pm 1}]]$ be an arbitrary Laurent series. Define
 {\em$\q$-derivation} of $a(z)$ as
 \begin{equation}\label{qder}
 \frac{d_{\q}}{d_{\q} z}a(z)=\frac{a(z\q)-a(z)}{z(\q-1)}\in
 V [[z^{\pm 1}]].
 \end{equation}
In order to simplify our notation we will denote the $n$th
$\q$-derivation of $a(z)$ as $a^{[n]}(z)$.

\begin{ex}
As an application of \eqref{qder} we calculate $\q$-derivation
of a monomial $a(z)=z^m$, $m\in\mathbb{Z}$:
\begin{equation*}
a^{[1]}(z)=(z^m)^{[1]}=\frac{(z\q)^m -z^m}{z(\q -1)}=\frac{\q^m
-1}{\q -1} z^{m-1}=[m]_{\q}z^{m-1}.
\end{equation*}
Specially, since $[0]_{\q}=0$, the $\q$-derivation of a constant equals zero.
\end{ex}

The operator $\frac{d_{\q}}{d_{\q}z}$ is obviously a linear operator.
Furthermore, it satisfies the general Leibniz rule:

\begin{pro}\label{qLeibniz}
For every nonnegative integer $m$ and $a(z),b(z)\in\om(V,V[[z]])$  
we have
\begin{equation}\label{qleibniz}
\left(a(z)b(z)\right)^{[m]}=\sum_{l=0}^{m}\gauss{m}{l}_{\q}a^{[l]}(z)b^{[m-l]}(z\q^l).
\end{equation}
\end{pro}

\begin{prf}
Since
$$\left(a(z)b(z)\right)^{[1]}=a(z)b^{[1]}(z)+a^{[1]}(z)b(z\q),$$
formula \eqref{qleibniz} follows by induction over $m$.
\end{prf}

\begin{rem}
It is easy to see that that \eqref{qleibniz} holds for any $a(z),b(z)\in\om(V,V((z)))$ satisfying 
$$a(z_1)b(z)\in \om(V,V((z_1,z))).$$
\end{rem}
Denote by $e^{z_0 \frac{\partial _{\q}}{\partial _{\q}z}}_{\q}$ the  operator
$$e^{z_0 \frac{\partial _{\q}}{\partial _{\q}z}}_{\q}=\sum_{l=0}^{\infty} \frac{1}{[l]_{\q}!} \frac{\partial_{\q}^l}{\partial_{\q} z^l} z_{0}^{l}\,\,\colon\, V[[z^{\pm 1}]] \to V[[z^{\pm 1}]][[z_0]],$$
where $z_0$  denotes right multiplication by the variable $z_0$ satisfying \eqref{wz}.
 At the end of this subsection we recall
$\q$-Taylor theorem:

\begin{pro}
For every $a(z)\in V[[z^{\pm 1}]]$ we have
\begin{equation}\label{qtaylor}
a(z+z_0) =e^{z_0 \frac{\partial _{\q}}{\partial _{\q}z}}_{\q}a(z),
\end{equation}
where $z_0$ and $z$ are subject to \eqref{wz}.
\end{pro}

\subsection{Quantum affine algebra \texorpdfstring{$U_{q}(\widehat{\mathfrak{sl}}_{2})$}{Uq(sl2)}}
For the simple Lie algebra $\mathfrak{sl}_2$, with the standard basis
$$x_{\alpha}=\begin{pmatrix}
0&1\\0&0\end{pmatrix},\quad x_{-\alpha}=\begin{pmatrix}
0&0\\1&0\end{pmatrix},\quad h_{\alpha}=\begin{pmatrix}
1&0\\0&-1\end{pmatrix},$$
denote by $\widehat{\mathfrak{sl}}_2$ the associated affine Lie algebra on the underlying vector space
$$\widehat{\mathfrak{sl}}_2=\mathfrak{sl}_2 \otimes \mathbb{C}[[t,t^{-1}]]\oplus\mathbb{C}c\oplus\mathbb{C}d,$$
with the bracket relations defined in a usual way (for details see \cite{Kac}).
For $a\in\mathfrak{sl}_2$ set 
$$a(z)=\sum_{r\in\mathbb{Z}} (a\otimes t^r)z^{-r-1}\in \widehat{\mathfrak{sl}}_2[[z,z^{-1}]].$$
Let $A=(a_{ij})_{i,j=0}^{1}$ be the generalized Cartan matrix   associated with 
$\widehat{\mathfrak{sl}}_{2}$ and
 $\widehat{\mathfrak{h}}\subset \widehat{\mathfrak{sl}}_{2}$  a vector space
over $\mathbb{C}$ with a basis
consisting of simple coroots $\alpha^{\vee}_{j}$, $j=0,1$, and
derivation $d$.
Denote by $\alpha_{0}, \alpha_{1}$ simple roots,
i.e. linear functionals from $\widehat{\mathfrak{h}}^{*}$  such that
$$\alpha_{i}(\alpha_{j}^{\vee})=a_{ji},\quad\alpha_{i}(d)=\delta_{i0},\quad
i,j=0,1,$$ 
and by $\Lambda_{0}, \Lambda_{1}$
fundamental weights, i.e.  elements of $\widehat{\mathfrak{h}}^{*}$ such that
$$\Lambda_{i}(\alpha_{j}^{\vee})=\delta_{ij},\quad
\Lambda_{i}(d)=0,\quad i,j=0,1.$$
Define a weight lattice $\widehat{P}$ of $\widehat{\mathfrak{sl}}_2$ as a free Abelian group generated by 
the elements $\Lambda_{0}, \Lambda_{1}$ and $\delta=\alpha_{0}+\alpha_{1}$.
An integral dominant weight is any   $\Lambda\in \widehat{P}$
such that  $\Lambda(\alpha_{i}^{\vee})\in\mathbb{Z}_{\geq 0}$  for
$i=0,1$. 
Denote by 
$$Q=\mathbb{Z}\alpha_{1}\subset\mathfrak{h}^{*}\quad\textrm{and}\quad P=\mathbb{Z}\lambda_{1}\subset\mathfrak{h}^{*}$$ 
the classical root lattice and the classical weight lattice of $\mathfrak{sl}_2$ respectively, where  $\lambda_{1}=\alpha_{1}/2$.

In this paper, we will mostly use $\q$-numbers introduced in the previous subsection.  However,
the definition of quantum affine algebra is usually given in terms of (differently defined) $q$-numbers, so we recall  them as well.
Fix an indeterminate $q$.
For any two integers $m$ and $l$, $l\geq 0$,   define $q$-integers,
 $$[m]=\frac{q^{m}-q^{-m}}{q-q^{-1}}, $$
  $q$-factorials,
 $$[0] !=1,\quad [l+1] !=[l+1][l]\cdots[1]$$
and $q$-binomial coefficients,
 $$\gauss{m}{l}=\frac{[m][m-1]\cdots [m-l+1]}{[l]!}.$$

We  recall  Drinfeld realization of the quantum affine algebra
$U_{q}(\widehat{\mathfrak{sl}}_{2})$.

\begin{defn}[\cite{D}]\label{drinfeld}
The quantum affine algebra $U_{q}(\widehat{\mathfrak{sl}}_{2})$ is the associative algebra over $\mathbb{C}(q^{1/2})$ with unit $1$ generated by the 
elements
$x^{\pm}(r)$, $a(s)$, $K^{\pm 1}$, $\gamma^{\pm 1/2}$ and $q^{\pm d}$,  $r,s\in\mathbb{Z}$, $s\neq 0$, subject to the following relations:
\begin{align}
& [\gamma^{\pm1/2},u]=0\textrm{ for all }u\in U_{q}(\widehat{\mathfrak{sl}}_{2}),\tag{D1}\label{D1}\\
& K K^{-1}=K^{-1}K=1,\tag{D2}\label{D2}\\
& [a(k),a(l)]=\delta_{k+l\hspace{2pt}0}\frac{[2k]}{k}\frac{\gamma^{k}-\gamma^{-k}}{q-q^{-1}},\tag{D3}\label{D3}\\
& [a(k),K^{\pm 1}]=[q^{\pm d},K^{\pm 1}]=0,\tag{D4}\label{D4}\\
& q^{d}x^{\pm}(k)q^{-d}=q^{k}x^{\pm}(k),\quad q^{d}a(k)q^{-d}=q^{k}a(k),\tag{D5}\label{D5}\\
& Kx^{\pm}(k)K^{-1}=q^{\pm 2}x^{\pm }(k) ,\tag{D6}\label{D6}\\
& [a(k),x^{\pm}(l)]=\pm\frac{[2k]}{k}\gamma^{\mp
|k|/2}x^{\pm}(k+l),\tag{D7}\label{D7}\\
& x^{\pm}(k+1)x^{\pm}(l)-q^{\pm 2}x^{\pm}(l)x^{\pm}(k+1)
=q^{\pm 2}x^{\pm}(k)x^{\pm}(l+1)-x^{\pm}(l+1)x^{\pm}(k),\tag{D8}\label{D8}\\
& [x^{+}(k),x^{-}(l)]=\frac{1}{q-q^{-1}}\left(\gamma^{\frac{k-l}{2}}\psi(k+l)-\gamma^{\frac{l-k}{2}}\phi(k+l)\right),\tag{D9}\label{D9}
\end{align}
where the elements $\phi(-r)$ and $\psi(r)$, $r\in\mathbb{Z}_{\geq 0}$, are given by 
\begin{align*}
& \phi(z)=\sum_{r=0}^{\infty}\phi(-r)z^{r}=K^{-1}\exp\left(-(q-q^{-1})\sum_{r=1}^{\infty}a(-r)z^{r}\right),\\
& \psi(z)=\sum_{r=0}^{\infty}\psi(r)z^{-r}=K\exp\left((q-q^{-1})\sum_{r=1}^{\infty}a(r)z^{-r}\right).
\end{align*}
\end{defn}

Denote by $x^{\pm}(z)$ the series
\begin{equation}\label{101_exp:series}
x^{\pm}(z)=\sum_{r\in\mathbb{Z}}x^{\pm}(r)z^{-r-1}\in U_{q}(\mathfrak{\widehat{\mathfrak{sl}}_{2}})[[z^{\pm 1}]].
\end{equation}
We shall continue to use the notation $x^{\pm}(z)$ for the action of  expression (\ref{101_exp:series}) on an arbitrary
$U_{q}(\widehat{\mathfrak{sl}}_{2})$-module $V$:
$$x^{\pm}(z)=\sum_{r\in\mathbb{Z}}x^{\pm}(r)z^{-r-1}\in (\ndo V)[[z^{\pm 1}]].$$


\subsection{Representations of \texorpdfstring{$U_{q}(\widehat{\mathfrak{sl}}_{2})$}{Uq(sl2)}}
First, we recall  Frenkel-Jing realization of  the integrable highest weight
 $U_{q}(\widehat{\mathfrak{sl}}_{2})$-modules $L(\Lambda_{i})$,
 $i=0,1$ (see \cite{FJ}).

Let $V$ be an arbitrary $U_{q}(\widehat{\mathfrak{sl}}_{2})$-module of
level $c$.
The Heisenberg algebra  $U_{q}(\widehat{\mathfrak{h}})$ of level $c$ is generated by the elements $a(k)$,  $k\in\mathbb{Z}\setminus\left\{0\right\}$, and the central element $\gamma^{\pm 1}=q^{\pm c}$ subject to the relations 
\begin{equation}\label{104_heisenberg}
[a(r),a(s)]=\delta_{r+s\hspace{2pt}0}\frac{[2r][cr]}{r},\quad
 r,s\in\mathbb{Z}\setminus\left\{0\right\}.
\end{equation}
Algebra $U_{q}(\widehat{\mathfrak{h}})$ has a natural realization on the space $\sym(\widehat{\mathfrak{h}}^{-})$
of the symmetric algebra generated by the elements $a(-r)$, 
  $r\in\mathbb{Z}_{>0}$, via the following rule:
  \begin{align*}
\gamma^{\pm 1}\hspace{5pt}&\ldots\hspace{5pt}\textrm{multiplication by }q^{\pm c},\\
a(r)\hspace{5pt}&\ldots\hspace{5pt}\textrm{differentiation operator subject to (\ref{104_heisenberg})},\\
a(-r)\hspace{5pt}&\ldots\hspace{5pt}\textrm{multiplication by the element }a_{i}(-r).
\end{align*}
Denote the resulted irreducible $U_{q}(\widehat{\mathfrak{h}})$-module
 as $K(c)$. Define the following operators on $K(c)$:
 \begin{align*}
 &E_{-}^{\pm}(a,z)=\exp\left(\mp\sum_{r\geq 1}\frac{q^{\mp cr/2}}{[cr]}a(-r)z^{r}\right),\\
  &E_{+}^{\pm}(a,z)=\exp\left(\pm\sum_{r\geq 1}\frac{q^{\mp cr/2}}{[cr]}a(r)z^{-r}\right).
\end{align*}

Let $\mathbb{C}\left\{Q\right\}$  be  group algebra of $Q$, which is generated by $e^{\alpha}$, $\alpha\in Q$. The space
$$\mathbb{C}\left\{P\right\}=\mathbb{C}\left\{Q\right\} \oplus \mathbb{C}\left\{Q\right\}e^{\lambda_1}$$
is a $\mathbb{C}\left\{Q\right\}$-module. Set 
\begin{equation}\label{remember}
L_0 =K(1)\otimes\mathbb{C}\left\{Q\right\},\qquad L_1 =K(1)\otimes\mathbb{C}\left\{Q\right\}e^{\lambda_1}.
\end{equation}
For $\alpha\in Q$ define an action  $z^{\partial_\alpha}$ on $\mathbb{C}\left\{P\right\}$ by
$$z^{\partial_\alpha}e^{\beta}=z^{(\alpha,\beta)}e^{\beta}.$$

\begin{thm}[\cite{FJ}]\label{frenkeljing}
 By the action 
\begin{align*}
x^{\pm}(z)&=E_{-}^{\pm}(-a,z)E_{+}^{\pm}(-a,z)\otimes e^{\pm\alpha}z^{\pm\partial_{\alpha}},
\end{align*}
the space $L_i$, $i=0,1$, becomes
the integrable highest weight module of $U_{q}(\widehat{\mathfrak{sl}}_{2})$  with the highest weight $\Lambda_i$.
\end{thm}

In \cite{Kozic} and \cite{Kozic2}, we considered some vertex operators on the space
$$L=K(1)\otimes \mathbb{C}\left\{P\right\},$$
 which were based on the explicit formulas for level $1$  vertex operators on integrable
highest weight modules of $U_{q}(\widehat{\mathfrak{sl}}_{n+1})$ found by Koyama in \cite{Koyama}. We now briefly recall Koyama's construction for $n=1$.
Define the following operators on the space $L$:
\begin{align*}
\mathcal{E}_{-}(z)&=\exp\left(\sum_{r=1}^{\infty}\frac{q^{r/2}}{[2r]}a(-r)z^{r}\right),\\
\mathcal{E}_{+}(z)&=\exp\left(-\sum_{r=1}^{\infty}\frac{q^{r/2}}{[2r]}a(r)z^{-r}\right).
\end{align*}
We define an operator 
$\mathcal{Y}(z)\in \om(L,L((z^{1/2})))$
by
\begin{align}\label{new-inter-}
\mathcal{Y}(z)=\mathcal{E}_{-}(z)\mathcal{E}_{+}(z)
\otimes e^{\lambda_1}(-z)^{\partial_{\lambda_1}}.
\end{align}
Applying the operator $\mathcal{Y}(z)$ on an arbitrary vector $v\in L$, we
get a formal power series in fractional powers $z^{1/2}$ of the
variable $z$, that has only a finite number of negative powers. Notice that  formula \eqref{new-inter-}  also defines an operator
 $L_0\to L_{1}((z))$.

The relations in the next proposition can be proved by a direct calculation. 

\begin{pro}\label{relations}
The following relations hold on $L=K(1)\otimes
\mathbb{C}\left\{P\right\}$:
{\allowdisplaybreaks\begin{align}
& x^{+}(z_1)x^{+}(z_2)=(z_1-z_2)(z_1-q^{-2}z_2):x^{+}(z_1)x^{+}(z_2):\label{r1}\\
& x^{+}(z_1)\mathcal{Y}(z_2)=\mathcal{Y}(z_2)x^{+}(z_1)=(z_1-z_2):x^{+}(z_1)\mathcal{Y}(z_2):\label{r2}\\
& x^{+}(z_1)\mathcal{E}_{-}(z_2)=\left(1-\frac{z_2}{z_1}\right)\mathcal{E}_{-}(z_2)x^{+}(z_1)\label{komutiranjesE-}\\
&E_{+}^{+}(-a,z_1)\mathcal{Y}(z_2)=\left(1-\frac{z_2}{z_1}\right)\mathcal{Y}(z_2)E_{+}^{+}(-a,z_1).\label{komutiranjesE2}
\end{align}}
\end{pro}

\section{Commutative operators and principal subspaces for \texorpdfstring{$U_{q}(\widehat{\mathfrak{sl}}_{2})$}{Uq(sl2)}}\label{komutativnisection}

\subsection{Commutative operators}\label{sec:1}

Define the
 following operators on the space
 $L$:
\begin{equation}\label{khat}\widehat{k}(z)=\exp\left((q-q^{-1})\sum_{r\geq
 1}\frac{-q^{r/2}}{1+q^{2r}}a(r)z^{-r}\right)\in \om(L,L ((z))).\end{equation}
Next,  define the operator  $\widehat{x}(z)$ by
\begin{equation}\label{new-op}
\widehat{x}(z)=x^{+}(z)\widehat{k}(z)\in \om(L,L ((z))).
\end{equation}
 
 \begin{rem}
Expression \eqref{new-op} is well-defined on an arbitrary
integrable $U_{q}(\widehat{\mathfrak{sl}}_{2})$-module. However, since the main results of
this section are proved for the level 1 case, we  consider
only the space $L$.
 \end{rem}

\begin{pro} On the space $L$ we have
\begin{equation}\label{komutativnost}
\widehat{x}(z_1)\widehat{x}(z_2)=\widehat{x}(z_2)\widehat{x}(z_1)=(z_1-z_2)^2
:\widehat{x}(z_1)\widehat{x}(z_2):.
\end{equation}
\end{pro}

Both of the above equalities can be proved by a direct calculation, using the
relations \eqref{D3}. Since the both sides of  \eqref{komutativnost} are
elements of $\om(L, L ((z_1,z_2)))$, we can apply the limit $\lim_{z_1 ,z_2
\to z}$ on \eqref{komutativnost}, thus getting

\begin{pro}
\begin{equation}\label{integrabilnost}
\widehat{x}(z)^2 = 0.
\end{equation}
\end{pro}

\begin{rem}
In \cite{DF}, certain vertex operators $\bar{x}^{\pm}_{i}(z)$, $i=1,2,...,n$, for the quantum affine algebra $U_q(\widehat{sl}_{n+1})$ were introduced. The level $1$ operators $\bar{x}^{+}(z)=\bar{x}^{+}_1 (z)$ for $U_q(\widehat{sl}_{2})$ are given by
$$\bar{x}^{+}(z)=x^{+}(z){k}^{+}(z)\in \om(L,L ((z))),$$
where
\begin{equation}\label{kplus}
k^{+}(z)=\exp\left((q-q^{-1})\sum_{r\geq 1}\frac{-q^{5r/2}}{1+q^{2r}}a(r)z^{-r}\right)\in \om(L,L ((z))).
\end{equation}
While the operators \eqref{khat} satisfy
$$(z_1 -q^{-2}z_2)\widehat{k}(z_1)x^{+}(z_2)=(z_1 -z_2)x^{+}(z_2)\widehat{k}(z_1),$$
the operators \eqref{kplus} satisfy
$$(z_1 -z_2)k^+ (z_1)x^{+}(z_2)=(z_1 -q^2z_2)x^{+}(z_2)k^+ (z_1).$$
As a result, the relations similar to \eqref{komutativnost}
hold for $\bar{x}^{+}(z)$ as well:
\begin{equation}\label{komutativnost2}
\bar{x}^+ (z_1)\bar{x}^+ (z_2)=\bar{x}^+ (z_2)\bar{x}^+ (z_1)=(z_1-q^{-2}z_2)(z_1-q^{2}z_2)
:\bar{x}^+(z_1)\bar{x}^+(z_2):.
\end{equation}
Since the both sides of  \eqref{komutativnost2} are
elements of $\om(L, L ((z_1,z_2)))$, we can apply the limit $\lim_{\substack{z_1 
\to z\\z_2 \to zq^{\pm 2}}}$ on \eqref{komutativnost2}, thus getting $q$-integrability relations
$$\bar{x}^+ (z)\bar{x}^+ (zq^{\pm 2})=0.$$
\end{rem}

Define the
 following operators on $L$: 
 \begin{align}
 &\widehat{\mathcal{E}}_{-}(z)=\exp\left(\sum_{r\geq
 1}\frac{q^{-r/2}}{2[r]}a(-r)z^r\right)\in \om(L,L ((z))),\nonumber\\
 & \widehat{\mathcal{E}}_{+}(z)=\exp\left(\sum_{r\geq
 1}\frac{-q^{r/2}}{[2r]}a(r)z^{-r}\right)\in\om(L,L ((z))),\nonumber\\
&\widehat{\mathcal{Y}}(z)=\widehat{\mathcal{E}}_{-}(z)\widehat{\mathcal{E}}_{+}(z)\otimes
e^{\lambda}(-z)^{\partial_\lambda}\in \om(L,L ((z^{1/2}))).\label{new-inter}
\end{align}
Notice that  formula \eqref{new-inter}  also defines an operator
 $L_0\to L_{1}((z))$.
Set 
$\widehat{\mathcal{E}}_{-,\lambda}(z)=\widehat{\mathcal{E}}_{-}(z)\otimes e^{\lambda}$.
 By a direct calculation one can prove
\begin{pro} 
On the space $L$ we have
\begin{align}
&\widehat{x}(z_1)\widehat{\mathcal{Y}}(z_2) = \widehat{\mathcal{Y}}(z_2)\widehat{x}(z_1)=(z_1 - z_2):\widehat{x}(z_1)\widehat{\mathcal{Y}}(z_2):,\label{dontneedthis}\\
&\widehat{x}(z_1)\widehat{\mathcal{E}}_{-,\lambda}(z_2)=-(z_1 - z_2)
\widehat{\mathcal{E}}_{-,\lambda}(z_2)\widehat{x}(z_1).\label{needthis}
\end{align}
\end{pro}
For a nonnegative integer $n$ set 
$$\widehat{x}^{(n)}(z)=\frac{\partial^n}{\partial
z^n}\widehat{x}(z)\qquad\text{and}\qquad
\widehat{\mathcal{E}}_{-,\lambda}^{(n)}(z)=\frac{\partial^n}{\partial
z^n}\widehat{\mathcal{E}}_{-,\lambda}(z).$$
By applying partial derivatives on \eqref{needthis}
we get
\begin{pro}
For nonnegative integers $n$ and $k$ we have
\begin{equation}\label{new-dif}
\widehat{x}^{(n)}(z_1)\widehat{\mathcal{E}}_{-,\lambda}^{(k)}(z_2)=k\widehat{\mathcal{E}}_{-,\lambda}^{(k-1)}(z_2)\widehat{x}^{(n)}(z_1)-n\widehat{\mathcal{E}}_{-,\lambda}^{(k)}(z_2)\widehat{x}^{(n-1)}(z_1)-(z_1-z_2
)\widehat{\mathcal{E}}_{-,\lambda}^{(k)}(z_2)\widehat{x}^{(n)}(z_1).
\end{equation}
\end{pro}

 \subsection{Basis for level 1 principal subspaces}\label{sec:2}

 Denote by $U_q (\widehat{\mathfrak{n}}^{\pm})$ a subalgebra of
 $U_{q}(\widehat{\mathfrak{sl}}_{2})$ generated by the elements $x^{\pm}(r)$,
 $r\in\mathbb{Z}$, and denote by $U_q (\widehat{\mathfrak{h}})_0$ a subalgebra
 of $U_{q}(\widehat{\mathfrak{sl}}_{2})$ generated by the elements $a(s)$,
 $K^{\pm 1}$, $\gamma^{\pm 1/2}$ and $q^{\pm d}$ for $s\in\mathbb{Z}$, $s\neq
 0$. It is well known that multiplication establishes an isomorphism of $\mathbb{C}(q^{1/2})$-vector
 spaces:
 $$U_{q}(\widehat{\mathfrak{sl}}_{2}) \cong U_q (\widehat{\mathfrak{n}}^{-})
 \otimes U_q (\widehat{\mathfrak{h}})_0 \otimes U_q
 (\widehat{\mathfrak{n}}^{+}).$$
 We recall the notion of  principal subspace from \cite{DF}. Let
 $\Lambda$ be an integral dominant weight. Denote by $v_{\Lambda}$ the highest
 weight vector of the integrable highest weight module $L(\Lambda)$. We define a
 principal subspace $W(\Lambda)$ of $L(\Lambda)$ as
 $$W(\Lambda)=U_q
 (\widehat{\mathfrak{n}}^{+}) v_{\Lambda}.$$
 In this subsection, we construct combinatorial bases for the principal
 subspaces $W(\Lambda_i)$, $i=0,1$, corresponding to the bases found in \cite{Kozic}. Each basis vector will be written as a
 monomial of endomorphisms $\widehat{x}(r)$ acting on the highest weight vector $v_{\Lambda_i}$.
 
Denote by  $\widehat{W}(\Lambda_i)$ a subspace of $L(\Lambda_i)$ spanned by the
vectors $$\widehat{x}(r_m)\cdots \widehat{x}(r_1)
v_{\Lambda_i},\quad\text{where}\quad m\in\mathbb{Z}_{\geq 0},\,r_j\in\mathbb{Z}.$$
The next Lemma can be proved in the same way as Lemma 9 in \cite{Kozic}.
\begin{lem} For $i=0,1$ 
$$\widehat{W}(\Lambda_i)=W(\Lambda_i).$$
\end{lem}

Equalities \eqref{komutativnost} and \eqref{integrabilnost} allow us to directly apply Georgiev's arguments (see \cite{G}). Hence, for $i=0,1$ we get the spanning set
$\widehat{\mathfrak{B}}_{W(\Lambda_i)}$ 
of the principal subspace $W(\Lambda_i)$, where
\begin{align}\label{baza}
\widehat{\mathfrak{B}}_{W(\Lambda_i)}=\bigg\{ \widehat{x}(l_m)\cdots
\widehat{x}(l_2)\widehat{x}(l_1)v_{\Lambda_i} \,\,\bigg| \,\,\,\bigg.
&l_1\leq -1-\delta_{1i}\text{ and }\bigg.
 l_{r+1}\leq l_r -2\\
\bigg. &\text{for all }l_r\in\mathbb{Z},\,r=1,2,\ldots,m-1,\,
m\in\mathbb{Z}_{\geq 0}\bigg\}.\nonumber
\end{align}
Since the difference and initial conditions in \eqref{baza}
coincide with the difference and initial conditions satisfied by  the elements of 
$\mathfrak{B}_{W(\Lambda_i)}\cdot v_{\Lambda_i}$, the basis of the principal subspace $W(\Lambda_i)$ found in \cite{Kozic},
we conclude
\begin{thm}
The set $\widehat{\mathfrak{B}}_{W(\Lambda_i)}$ forms a basis for
the principal subspace $W(\Lambda_i)$.
\end{thm}

\subsection{Vertex algebra \texorpdfstring{$W_0$}{W0}}\label{sec:3}

Let $W_0\subset\om(L_0,L_0 ((z)))$ be the vertex algebra generated by
$\widehat{x}(z)\in\om(L_0,L_0 ((z)))$.
 Denote by
$$Y(\cdot,z_0)\colon W_0\to (\ndo W_0)[[z_0,z_{0}^{-1}]]$$
the corresponding vertex operator map.
Details about construction theorems for vertex algebras and their modules can be
found in \cite{LiLep}.
For every $a(z),b(z)\in W_0$ we have
$$Y(a(z),z_0)b(z)=\sum_{r\in\mathbb{Z}} a(z)_r b(z)z_{0}^{-r-1}=a(z+z_0)b(z)$$
and, therefore,
\begin{equation}\label{classical}
a(z)_r b(z) =\left\{\begin{array}{l@{\,\ }l}
0& \text{if }r\geq 0,\\
\frac{1}{(-r-1)!}a^{(-r-1)}(z)b(z)&\text{if }r<0.\end{array}\right.
\end{equation}

\begin{rem}
Since the operator $\widehat{x}(z)$ is commutative, we can consider the smallest commutative, associative algebra $A_0$ with unit $1$ and derivation $\frac{d}{dz}$,
which contains  $\widehat{x}(z)$. The algebra $A_0$ obviously coincides with $W_0$ (as a vector space). Furthermore, the vertex algebra structure on $W_0$ can be recovered from $A_0$
by Borcherds' construction (see \cite{Bor}).
\end{rem}

The space $W_0$ is spanned by the vectors
\begin{equation}\label{span1}
\widehat{x}(z)_{l_m}\ldots\widehat{x}(z)_{l_2}\widehat{x}(z)_{l_1}1,\quad
l_m\leq\ldots\leq l_2\leq l_1\leq -1,\,\,m\in\mathbb{Z}_{\geq 0}.
\end{equation}
By applying Leibniz rule on
$$\frac{d^{n}}{dz^n}\left(\widehat{x}(z)\widehat{x}(z)\right)=\frac{d^{n}}{dz^n} 0=0$$
we get
\begin{equation}\label{derivation765}
\sum_{l=0}^{n}\binom{n}{l}\widehat{x}^{(l)}(z)\widehat{x}^{(n-l)}(z)=0.
\end{equation}
This equality, together with \eqref{komutativnost}, allows us to express a 
product
$$\widehat{x}(z)_s \widehat{x}(z)_r,\quad s-r=0,\pm 1,$$
as a linear combination of products
$$\widehat{x}(z)_t \widehat{x}(z)_u,\quad u-t\geq 2,\,\, u+t=s+r.$$
Therefore, we can reduce  spanning set  \eqref{span1},
thus getting a smaller spanning set of $W_0$:
\begin{align}\label{W0baza}
\widehat{\mathfrak{B}}_{W_0}=\bigg\{ \widehat{x}(z)_{l_m}\ldots\widehat{x}(z)_{l_2}\widehat{x}(z)_{l_1}1 \,\,\bigg| \,\,\,\bigg.
&l_1\leq -1\text{ and }\bigg.
 l_{r+1}\leq l_r -2\\
\bigg. &\text{for all }l_r\in\mathbb{Z},\,r=1,2,\ldots,m-1,\,
m\in\mathbb{Z}_{\geq 0}\bigg\}.\nonumber
\end{align}
\begin{rem}
Notice that equality \eqref{komutativnost}  implies
\begin{equation}\label{niceexample1}
\widehat{x}(z)_{-1}\widehat{x}(z)_{-1}1=\widehat{x}(z)_{-2}\widehat{x}(z)_{-1}1=\widehat{x}(z)_{-2}\widehat{x}(z)_{-1}1=0,
\end{equation}
while \eqref{derivation765} can be written as
\begin{equation}\label{niceexample2}
n! \sum_{l=0}^{n} \widehat{x}(z)_{-l-1} \widehat{x}(z)_{-n+l-1} 1=0.
\end{equation}
\end{rem}

Our next goal is to prove the linear independence of $\widehat{\mathfrak{B}}_{W_0}$.

\begin{lem}\label{zeroth}
For every positive integer $m$ and $\varepsilon =0,1$
$$\widehat{x}(z)_{-2m+1-\varepsilon}\widehat{x}(z)_{-2m+1}\ldots\widehat{x}(z)_{-3}\widehat{x}(z)_{-1}1=
0.$$
\end{lem}

\begin{prf}
The Lemma is a consequence of equalities \eqref{komutativnost} and it can be proved by induction on $m$.
\end{prf}

\begin{lem}\label{first}
For every positive integer $m$
$$\widehat{x}(z)_{-2m+1}\ldots\widehat{x}(z)_{-3}\widehat{x}(z)_{-1}1\neq 0.$$
\end{lem}

\begin{prf} Set
$$
c_0 =\frac{1}{(2m-2)!\cdots 4! 2!},\qquad
d_0=\frac{\partial ^{2m-2}}{\partial z_{m}^{2m-2}}\cdots
\frac{\partial ^{2}}{\partial z_{2}^{2}}\left(\prod_{\substack{j,k=1\\j>k}}^{m}(z_j
- z_k)^2\right)\in\mathbb{Z}\setminus\left\{0\right\}.
$$
We have
\begin{align*}
&\widehat{x}(z)_{-2m+1}\ldots\widehat{x}(z)_{-3}\widehat{x}(z)_{-1}1=c_0 \widehat{x}^{(2m-2)}(z)\ldots\widehat{x}^{(2)}(z)\widehat{x}(z)\\
&=c_0\lim_{z_m,\ldots,z_2 , z_1\to
z}\left(\widehat{x}^{(2m-2)}(z_m)\cdots\widehat{x}^{(2)}(z_2)\widehat{x}(z_1)\right)\\
&=c_0\lim_{z_m,\ldots,z_2 , z_1\to
z}\left(\frac{\partial ^{2m-2}}{\partial z_{m}^{2m-2}}\cdots
\frac{\partial ^{2}}{\partial z_{2}^{2}}\left(\widehat{x}(z_m)\cdots\widehat{x}(z_2)\widehat{x}(z_1)\right)\right)\\
&=c_0\lim_{z_m,\ldots,z_2 , z_1\to
z}\left(\frac{\partial ^{2m-2}}{\partial z_{m}^{2m-2}}\cdots
\frac{\partial ^{2}}{\partial z_{2}^{2}}\left(\left(\prod_{\substack{j,k=1\\j>k}}^{m}(z_j
- z_k)^2\right)
:\widehat{x}(z_m)\ldots\widehat{x}(z_2)\widehat{x}(z_1):\right)\right)\\
&=c_0\lim_{z_m,\ldots,z_2 , z_1\to
z}\left(\left(\frac{\partial ^{2m-2}}{\partial z_{m}^{2m-2}}\cdots
\frac{\partial ^{2}}{\partial z_{2}^{2}}\left(\prod_{\substack{j,k=1\\j>k}}^{m}(z_j
- z_k)^2\right)\right)
:\widehat{x}(z_m)\ldots\widehat{x}(z_2)\widehat{x}(z_1):\right)\\
&=c_0 d_0 \lim_{z_m,\ldots,z_2 , z_1\to
z}
:\widehat{x}(z_m)\ldots\widehat{x}(z_2)\widehat{x}(z_1):\\
&= c_0 d_0 :\underbrace{\widehat{x}(z)\ldots\widehat{x}(z)\widehat{x}(z)}_{m}:\,\neq\, 0.
\end{align*}
\end{prf}

\begin{lem}\label{second}
The elements of the set $\widehat{\mathfrak{B}}_{W_0}$ are nonzero.
\end{lem}

\begin{prf}
Suppose that there exists
\begin{equation}\label{help}
b(z)=\widehat{x}(z)_{l_m}\ldots\widehat{x}(z)_{l_2}\widehat{x}(z)_{l_1}1\in
\widehat{\mathfrak{B}}_{W_0}
\end{equation} such that $b(z)=0$, i.e.
\begin{equation}\label{31415}
\widehat{x}^{(-l_m-1)}(z)\cdots\widehat{x}^{(-l_2-1)}(z)\widehat{x}^{(-l_1
-1)}(z)=0.
\end{equation}
We can assume that 
\begin{equation}\label{help2}
(-l_m-1,\ldots,-l_2 -1,-l_1-1)
\end{equation} is a minimal $m$-tuple
in the lexicographic ordering such that
$\widehat{x}(z)_{l_m}\ldots\widehat{x}(z)_{l_2}\widehat{x}(z)_{l_1}1=0$, so
if
$$c(z)=\widehat{x}(z)_{k_m}\ldots\widehat{x}(z)_{k_2}\widehat{x}(z)_{k_1}1=0$$
 and $c(z) \in \widehat{\mathfrak{B}}_{W_0}$, then $$(-k_m-1,\ldots,-k_2
-1,-k_1-1)\geq (-l_m-1,\ldots,-l_2 -1,-l_1-1).$$
Notice that Lemma \ref{first} implies
$$(-l_m-1,\ldots,-l_2 -1,-l_1-1)>(2m-2,\ldots,2,0).$$
Therefore, there exists an integer $k\in\left\{1,2,\ldots,m\right\}$ such that
$-l_k-1>2(k-1)$. Let $j\in\left\{1,2,\ldots,m\right\}$ be the minimal integer such
that $-l_j-1>2(j-1)$.
Equality \eqref{31415} implies
\begin{equation*}
\widehat{\mathcal{E}}_{-,\lambda}(z_1)^{-1}\widehat{x}^{(-l_m-1)}(z)\cdots\widehat{x}^{(-l_2-1)}(z)\widehat{x}^{(-l_1
-1)}(z)\widehat{\mathcal{E}}_{-,\lambda}^{(j-1)}(z_1)=0,
\end{equation*}
where
$$\widehat{\mathcal{E}}_{-,\lambda}(z_1)^{-1}=\widehat{\mathcal{E}}_{-}(z_1)^{-1}\otimes e^{-\lambda}.$$
Therefore,
\begin{align*}
0=\lim_{z_1 \to z} \widehat{\mathcal{E}}_{-,\lambda}(z_1)^{-1}\widehat{x}^{(-l_m-1)}(z)\cdots\widehat{x}^{(-l_2-1)}(z)\widehat{x}^{(-l_1
-1)}(z)\widehat{\mathcal{E}}_{-,\lambda}^{(j-1)}(z_1).
\end{align*}
Now, we use \eqref{new-dif} and Lemma \ref{zeroth} in order to move the operator $\widehat{\mathcal{E}}_{-,\lambda}^{(j-1)}(z_1)$ to the left:
\begin{align*}
0=\lim_{z_1 \to z} (j-1)! \widehat{\mathcal{E}}_{-,\lambda}(z_1)^{-1}\widehat{x}^{(-l_m-1)}(z)\cdots\widehat{x}^{(-l_j-1)}(z)
\widehat{\mathcal{E}}_{-,\lambda}(z_1)\widehat{x}^{(-l_{j-1}-1)}(z)\cdots
\widehat{x}^{(-l_1-1)}(z).
\end{align*}
Finally, we have
\begin{align*}
0&=\lim_{z_1 \to z} c_0 \widehat{\mathcal{E}}_{-,\lambda}(z_1)^{-1}\widehat{\mathcal{E}}_{-,\lambda}(z_1)\widehat{x}^{(-l_m-2)}(z)\cdots\widehat{x}^{(-l_j-2)}(z)
\widehat{x}^{(-l_{j-1}-1)}(z)\cdots
\widehat{x}^{(-l_1-1)}(z)\\
&=\lim_{z_1 \to z} c_0\widehat{x}^{(-l_m-2)}(z)\cdots\widehat{x}^{(-l_j-2)}(z)
\widehat{x}^{(-l_{j-1}-1)}(z)\cdots
\widehat{x}^{(-l_1-1)}(z)\\
&=c_0\widehat{x}^{(-l_m-2)}(z)\cdots\widehat{x}^{(-l_j-2)}(z)
\widehat{x}^{(-l_{j-1}-1)}(z)\cdots
\widehat{x}^{(-l_1-1)}(z),
\end{align*}
where
$$c_0=(-l_m-1)\cdots (-l_{j+1} -1)(-l_j
-1)\cdot (j-1)!\neq 0.$$
Hence, we have
$$\widehat{x}^{(-l_m-2)}(z)\cdots\widehat{x}^{(-l_j-2)}(z)\widehat{x}(z)^{(-l_{j-1}-1)}\cdots\widehat{x}(z)^{(-l_2-1)}\widehat{x}^{(-l_1
-1)}(z)=0.$$
Since
$$\widehat{x}(z)_{l_m
+1}\ldots\widehat{x}(z)_{l_j
+1}\widehat{x}(z)_{l_{j-1}
}\ldots\widehat{x}(z)_{l_2}\widehat{x}(z)_{l_1}1\in
\widehat{\mathfrak{B}}_{W_0}$$
and
\begin{align*}
&(-l_m-1,\ldots,-l_j -1, -l_{j-1},\ldots,-l_2 -1,-l_1-1)\\
&\qquad\qquad >(-l_m-2,\ldots,-l_j -2, -l_{j-1}-1,\ldots,-l_2
-1,-l_1-1),
\end{align*}
we have a contradiction to a minimality of $j$.
\end{prf}

Now, we can prove the main result of this section.

\begin{thm}\label{findthis}
The set $\widehat{\mathfrak{B}}_{W_0}$ forms a basis for $W_0$.
\end{thm}

\begin{prf}
Consider the linear combination 
\begin{equation}\label{linear0}
\sum_{j=1}^{n}\alpha_{j}b_{j}(z)=0,
\end{equation}
where
$\alpha_j$ are nonzero scalars and
$$b_{j}(z)=\widehat{x}(z)_{l_{m_j}^{(j)}}\ldots\widehat{x}(z)_{l_2^{(j)}}\widehat{x}(z)_{l_1^{(j)}}1\in
\widehat{\mathfrak{B}}_{W_0}.$$
We can assume that $n$ is the smallest positive integer for which such a linear
combination exists. Lemma \ref{second}  implies that $n\geq 2$. 
 Furthermore, we can assume that $m_j=m_k$ for all $j,k=1,2,...,n$. Indeed, if $m_j\neq m_k$ 
   for some $j,k=1,2,...,n$, by applying the operator $1\otimes q^\alpha$ on \eqref{linear0} we get
     $$\sum_{j=1}^{n} q^{2m_j}\alpha_j b_{j}(x) =0,$$
     which, together with \eqref{linear0}, contradicts to minimality of $n$.
 Therefore, without
loss of generality we can assume that $m_j=m_k$ for all
$j,k=1,2,\ldots,n$.

In the proof of Lemma \ref{second} we described a procedure of reducing an
element
$$\widehat{x}(z)_{l_m}\ldots\widehat{x}(z)_{l_j
}\widehat{x}(z)_{l_{j-1}
}\ldots\ldots\widehat{x}(z)_{l_2}\widehat{x}(z)_{l_1}1\in
\widehat{\mathfrak{B}}_{W_0}$$ up to an element
$$\widehat{x}(z)_{l_m
+1}\ldots\widehat{x}(z)_{l_j
+1}\widehat{x}(z)_{l_{j-1}
}\ldots\widehat{x}(z)_{l_2}\widehat{x}(z)_{l_1}1\in
\widehat{\mathfrak{B}}_{W_0}.$$
We can continue to apply such a procedure until we get the ''minimal element"
\begin{equation}\label{maxisss}\widehat{x}(z)_{-2m+1}\ldots\widehat{x}(z)_{-2j+1
}\widehat{x}(z)_{-2(j-1)+1}
\ldots\widehat{x}(z)_{-3}\widehat{x}(z)_{-1}1\in
\widehat{\mathfrak{B}}_{W_0}.
\end{equation}
In the same proof  we also associated with  element
\eqref{help} of $\widehat{\mathfrak{B}}_{W_0}$  $m$-tuple \eqref{help2}. For two such
elements $b(z),c(z)\in \widehat{\mathfrak{B}}_{W_0}$ we shall write $b(z)<c(z)$
if such an inequality holds for their corresponding $m$-tuples. 

Suppose that for the operators $b_{j}(z)$ in \eqref{linear0} we have
$b_{n}(z)<\ldots<b_{1}(z)$. We can keep applying our reduction procedure
on \eqref{linear0}, not stopping until $b_{1}(z)$ does not become  minimal
(nonzero) element  \eqref{maxisss}. Notice that all the other
elements $b_{j}(z)$, $j>1$, were already annihilated at some intermediate stage, so we get
$$c_0\alpha_1 b_1 (z)=0$$
for some nonzero scalar $c_0$. This implies $\alpha_1=0$. Contradiction.
\end{prf}

The space
$$W_1 =\left\{a(z)\widehat{\mathcal{Y}}(z)\,\,\big|\big.\,
a(z)\in W_0\right\} \subset\om(L_0, L_{1}((z)))$$
can be equipped with a $W_0$-module structure by choosing a linear map
$$Y_{W_1}(\cdot,z_0)\colon W_0\to \left(\ndo W_1 \right)  [[z^{\pm 1}_{0}]]$$
such that
$$Y_{W_1}(a(z),z_0)b(z)=a(z+z_0)b(z)$$
for all $a(z)\in W_0$, $b(z)\in W_1$.
Set
\begin{align}\label{W1baza}
\widehat{\mathfrak{B}}_{W_1}=\bigg\{ \widehat{x}(z)_{l_m}\ldots\widehat{x}(z)_{l_2}\widehat{x}(z)_{l_1}\widehat{\mathcal{Y}}(z) \,\,\bigg| \,\,\,\bigg.
&l_1\leq -2\text{ and }\bigg.
 l_{r+1}\leq l_r -2\\
\bigg. &\text{for all }l_r\in\mathbb{Z},\,r=1,2,\ldots,m-1,\,
m\in\mathbb{Z}_{\geq 0}\bigg\}.\nonumber
\end{align}
Condition $l_1\leq -2$ is a consequence of  \eqref{dontneedthis} because by applying the limit $\lim_{z_2 ,z_1 \to z}$ on this relation
 we get
 $$\widehat{x}(z)_{-1}\widehat{\mathcal{Y}}(z)=\widehat{x}(z)\widehat{\mathcal{Y}}(z)=0.$$

\begin{thm}
The set $\widehat{\mathfrak{B}}_{W_1}$ forms a basis for $W_1$.
\end{thm}

\begin{prf}
The proof is similar to the proof of Theorem \ref{findthis}, so we only briefly outline it.
First, $\widehat{\mathfrak{B}}_{W_1}$ obviously spans $W_1$.
Next, relations \eqref{komutativnost} and \eqref{dontneedthis} imply
$$\widehat{x}(z_m)\cdots \widehat{x}(z_1)\widehat{\mathcal{Y}}(z) =\prod_{\substack{j,k=1\\j>k}}^{m}(z_j
- z_k)^2 \prod_{\substack{j=1}}^{m}(z_j
- z):\widehat{x}(z_m)\ldots \widehat{x}(z_1)\widehat{\mathcal{Y}}(z):,$$
so we can proceed similarly as in the proofs of Lemmas \ref{first}, \ref{second} and Theorem \ref{findthis}.
However, since
$$\widehat{\mathcal{Y}}(z)\widehat{\mathcal{E}}_{-,\lambda}(z_1)=\left(1-\frac{z_1}{z}\right)^{1/2}z^{1/2}\widehat{\mathcal{E}}_{-,\lambda}(z_1)\widehat{\mathcal{Y}}(z),$$
the operator $\widehat{\mathcal{E}}_{-,\lambda}^{(j-1)}(z_1)$, which was used in the proof of Lemma \ref{second}, should be replaced by
 $$\frac{\partial^{j-1}}{\partial z_{1}^{j-1}}\left(\left(1-\frac{z_1}{z}\right)^{-1/2}z^{-1/2}\widehat{\mathcal{E}}_{-,\lambda}(z_1)\right).$$
 \end{prf}

\section{Nonlocal \texorpdfstring{$\q$}{q}-vertex algebras}\label{teorijskisection}

\subsection{Definition}

Let $L$ be an arbitrary vector space over the field $\mathbb{F}\supseteq \mathbb{C}(\q)$ of characteristic zero. Denote by $1=1_L$ the identity
$L\to L$. In this section we consider  certain vertex algebra-like structures on subspaces of
$$\mathcal{E}(L)=\om(L,L((z))).$$

\begin{defn}\label{kva}
A {\em nonlocal $\q$-vertex algebra} is a 3-tuple $(V,Y,1)$, where  $V\subseteq \mathcal{E}(L)$ is a
vector space equipped with a linear map
\begin{align}
Y(\cdot,z_0)\colon V &\to (\ndo V)[[z^{\pm 1}_{0}]]\nonumber\\
a(z)&\mapsto Y(a(z),z_0 )=\sum_{r\in\mathbb{Z}}a(z)_{r}z^{-r-1}_{0}\tag{v0}\label{V0}
\end{align}
\noindent and  with a distinguished vector $1\in V$ such that the
following conditions hold:
\noindent For every $a(z),b(z),c(z)\in V$
\begin{align}
&a(z)_r b(z)=0\text{ for }r\text{ sufficiently large;}\tag{v1}\label{V1}\\
&Y(1,z_0)=1_{V};\tag{v2}\label{V2}\\
&Y(a(z),z_0)1\in V[[z_0]]\quad\text{and}\quad \lim_{z_0\to
0}Y(a(z),z_0)1=a(z);\tag{v3}\label{V3}\\
&Y(a(z),z_0 +z_2)Y(b(z),z_2) c(z)=Y(Y(a(z),z_0)b(z),z_2)c(z),\tag{v4}\label{V4}
\end{align}
where
\begin{equation}\tag{v5}\label{qvar}
z_2 z_0=\q z_0 z_2,\quad z_2 z=\q zz_2,\quad z_0 z=\q z_0 z.
\end{equation}
\end{defn}
The definition requires some further explanations.
In \eqref{V0}, by $z_0$ is denoted right multiplication by $z_0$, i.e. for $a(z),b(z) \in V$ we have
$$ Y(a(z),z_0 )b(z)=\sum_{r\in\mathbb{Z}}a(z)_{r} b(z)z^{-r-1}_{0}\in V((z_{0})).$$

Since the variables $z_2,z_0,z$ do not commute, we use the following conventions in  \eqref{V4}. On the left-hand side of the equality
we assume that $z_0 +z_2$ appears to the right of $z_2$:
\begin{align*}
 Y(a(z),z_0 +z_2)Y(b(z),z_2) c(z)=\sum_{r\in\mathbb{Z}}\sum_{s\in\mathbb{Z}} a(z)_r ( b(z)_s c(z)) z_2^{-s-1} (z_0 +z_2)^{-r-1}.
\end{align*}
Of course, we expand $(z_0 +z_2)^{-r-1}$ in nonnegative powers of $z_2$ as in  Proposition \ref{expand}.

When expanding  the right-hand side of \eqref{V4}, we  assume that variable $z_0$ appears to the left of $c(z)$:
\begin{align*}
&Y(Y(a(z),z_0)b(z),z_2)c(z)=Y(\sum_{r\in\mathbb{Z}}a(z)_r b(z)z_0^{-r-1},z_2)c(z)\\
&\qquad=\sum_{s\in\mathbb{Z}} \sum_{r\in\mathbb{Z}}(a(z)_r b(z))_s z_0^{-r-1}c(z) z_2^{-s-1}\\
&\qquad=\sum_{s\in\mathbb{Z}} \sum_{r\in\mathbb{Z}}(a(z)_r b(z))_s c(z\q^{-r-1}) z_0^{-r-1} z_2^{-s-1}.
\end{align*}

\begin{defn}\label{kvam}
Let $(V,Y,1)$ be a nonlocal $\q$-vertex algebra and $K_1,K_2$ vector spaces over the field $\mathbb{F}$. A $V${\em-module} is a vector space $W\subset\om(K_1 ,K_2 ((z)))$ equipped with  linear map $Y_W\colon
V\to\om(W,W((z_0)))$  such that for every $a(z),b(z)\in V$, $c(z)\in W$ we have
\begin{align}
& Y_{W}(1,z_0)=1_{W};\tag{m1}\label{M1}\\
& Y_W (a(z),z_0 +z_2)Y_W(b(z),z_2)c(z)=Y_W
(Y(a(z),z_0)b(z),z_2)c(z),\tag{m2}\label{M2}
\end{align}
where $z_2 ,z_0,z$ satisfy 
\eqref{qvar}.
\end{defn}

Relation \eqref{M2} should be understood in the same way as \eqref{V4}.

\subsection{Construction of nonlocal \texorpdfstring{$\q$}{q}-vertex algebras}

The main goal of this subsection is providing a construction method of nonlocal $\q$-vertex
algebras and their modules. Even though the main result, together with its proof, is motivated by the results in \cite{LiNonlocal},
we are studying  some different products among operators, which will prove to be useful in dealing with principal subspaces.
Unless stated otherwise, we assume that all the variables $z,z_0,z_1,z_2,...$ are commutative. 
 Denote by $R = R_{\q}$ the operator $\mathcal{E}(L)\to \mathcal{E}(L)$ given by
$$R a(z)=a(z\q).$$
We shall consider  vertex operators which satisfy the following  special case of quasi-compatibility:

\begin{defn}\label{qkomm}
A sequence  $(a_1 (z),\ldots,a_m (z))$ in
$\mathcal{E}(L)$ is said to be {\em quasi-commutative} if
\begin{equation}\label{qcomm}
a_1(z_1)a_2 (z_2)\cdots a_m(z_m)\in\om(L,L((z_1,z_2,\ldots,z_m))).
\end{equation}
\end{defn}

Let $a(z),b(z)$ be operators in $\mathcal{E}(L)$ such that $$[a(z_1),b(z)]=0.$$
Then the pairs $(a(z),b(z))$ and $(b(z),a(z))$ are 
quasi-commutative.

\begin{defn}\label{produkti}
Let $(a(z),b(z))$ be a quasi-commutative pair in $\mathcal{E}(L)$. For
$r\in\mathbb{Z}$ we define $a(z)_r b(z)\in (\ndo L)[[z^{\pm 1}]]$ by
\begin{equation}\label{rtiprodukti}
 a(z)_r b(z)=\left\{\begin{array}{l@{\,\ }l}
0& \text{if }r\geq 0,\\
\frac{1}{[-r-1]_{\q}!}a^{[-r-1]}(z)b(z\q^{-r-1})&\text{if }r<0.\end{array}\right.
\end{equation}
\end{defn}

Recall \eqref{classical} and notice that, for commutative operators, the  limit $\q\to 1$ of the
$r$th products given by Definition \ref{produkti} coincides with the standard
$r$th products of local vertex operators.

First, we notice that the space $\mathcal{E}(L)$ is closed under the operations
given by the preceding definition:

\begin{pro}
Let $(a(z),b(z))$ be a quasi-commutative pair in $\mathcal{E}(L)$.
 Then
$$a(z)_r b(z)\in \mathcal{E}(L)\quad\text{for every integer }r.$$
\end{pro}
Next, we list  some elementary properties of the $r$th products: 

\begin{pro}\label{bigprop}
Let $(a(z),b(z))$ be a quasi-commutative pair in $\mathcal{E}(L)$.
\begin{enumerate}[(a)]
  \item For every $\alpha,\beta\in\mathbb{F}$ the pair
  $(a(z\alpha),b(z\beta))$ is quasi-commutative.
  \item For every  $r,s\geq 0$ the pair $(a^{[r]}(z),b^{[s]}(z))$ is
  quasi-commutative.
  \item Pair $(1_L,a(z))$ is quasi-commutative and for every integer $r$ we
  have \begin{equation}\label{vv2}
  (1_L)_r a(z)=\delta_{r+1\, 0} a(z).
  \end{equation}
  \item Pair $(a(z),1_L)$ is quasi-commutative and we have 
  \begin{equation}\label{vv3}
  a(z)_{-1} 1_L=a(z).
  \end{equation}
\end{enumerate}
\end{pro}

For the operator $a(z)\in\mathcal{E}(L)$ we  write
\begin{align}\label{Y}
Y(a(z),z_0)&=\sum_{r=0}^{\infty} a(z)_{-r-1} z_{0}^{r}
=\sum_{r=0}^{\infty}\frac{1}{[r]_{\q}!}a^{[r]}(z) z_{0}^{r},
\end{align}
where, as before, by $z_0$ we denote the right multiplication with variable $z_0$ satisfying
\begin{equation}\label{asbefore}
z_0 z = \q zz_0.
\end{equation}
By applying \eqref{Y} on an  operator $b(z)\in\mathcal{E}(L)$, such that the pair $(a(z),b(z))$ is quasi-commutative, we get
\begin{align}\label{YY}
Y(a(z),z_0)b(z)&=\sum_{r=0}^{\infty} a(z)_{-r-1} b(z) z_{0}^{r}
=\sum_{r=0}^{\infty}\frac{1}{[r]_{\q}!}a^{[r]}(z) b(z\q^r ) z_{0}^{r}\in\mathcal{E}(L)[[z_0]].
\end{align}

\begin{pro}\label{vvv23}
For every quasi-commutative pair $(a(z),b(z))$ in $\mathcal{E}(L)$ we have
\begin{align}
&Y(1_L,z_0)=1_{\mathcal{E}(L)};\label{vvv2}\\
&Y(a(z),z_0)1_{L}\in\mathcal{E}(L)[[z_0]]\quad\text{and}\quad\lim_{z_0\to
0}Y(a(z),z_0)1_{L}=a(z);\label{vvv3}\\
&Y(a^{[1]}(z),z_0)b(z\q)=\frac{\partial_{\q}}{\partial_{\q}
z_0}\left(Y(a(z),z_0)b(z )\right).\label{nolabel}
\end{align}
\end{pro}
\begin{prf}
Equalities \eqref{vvv2} and \eqref{vvv3} follow from \eqref{rtiprodukti}, \eqref{vv2}
 and \eqref{vv3}, while equality \eqref{nolabel} can be easily
verified by a direct calculation.
\end{prf}

Using statement (a) of Proposition \ref{bigprop} one can prove:

\begin{lem}\label{dong}
Let  $(a_1(z),a_2(z), ... ,a_m(z))$ and $(a_{j}(z),a_{j+1}(z))$, $j=1,2,...,m-1$, be  quasi-commutative sequences in $\mathcal{E}(L)$. Then, for every integer $r$    the sequence
$$(a_1(z),...,a_{j-1}(z) ,a_{j}(z)_r a_{j+1}(z)  ,a_{j+2}(z),...  ,a_{m}(z)  )$$ is quasi-commutative.
\end{lem}

\begin{rem}
The operator $a_{j}(z)_r a_{j+1}(z)$ can be expressed as a $\mathbb{F}[z^{\pm 1}]$-linear combination of certain minus first products, introduced in \cite{LiNonlocal},  of the operators $a_j (z\q^{m}),m\in\mathbb{Z}$, and the operator $a_{j+1}(z)$. Therefore, Lemma \ref{dong} is a consequence of    \cite[Proposition 2.18.]{LiNonlocal}.
\end{rem}

Let $\mathcal{S}$ be a subset of $\mathcal{E}(L)$. We shall say that $\mathcal{S}$ is {\em
quasi-commutative} if every finite sequence in $\mathcal{S}$ is
quasi-commutative. We shall say that $\mathcal{S}$ is {\em $closed$} if
\begin{align*}
& \,a(z)_r b(z)\in \mathcal{S}\quad\text{for all }a(z),b(z)\in \mathcal{S},\, r\in\mathbb{Z}.
\end{align*}

\begin{thm}\label{construction1}
Let $V$ be a closed and quasi-commutative subspace of $\mathcal{E}(L)$ such
that $1_{L}\in V$. Then $(V,Y,1_{L})$
is a nonlocal $\q$-vertex algebra.
\end{thm}

\begin{prf}
Assertion \eqref{V1} follows from
\eqref{YY}, while \eqref{V2} and \eqref{V3} have already been proved in Proposition
\ref{vvv23}. Suppose that variables $z_2,z_0,z$ satisfy \eqref{qvar}.
For $a(z),b(z),c(z)\in V$,  using 
\eqref{qbinomial} and
\eqref{qleibniz}, we get
\begin{align}
&Y(Y(a(z),z_0)b(z),z_2)c(z)
=Y(\sum_{r=
0}^{\infty}\frac{1}{[r]_{\q}!}a^{[r]}(z)b(z\q^r )z_0^{r},z_2)c(z)\nonumber\\
&\qquad=\sum_{s= 0}^{\infty}\sum_{r= 0}^{\infty}
\frac{1}{[r]_{\q}![s]_{\q}!}(a^{[r]}(z)b(z\q^r))^{[s]}z_0^{r}c(z\q^s )z_2^{s}\nonumber\\
&\qquad=\sum_{s= 0}^{\infty}\sum_{r=
0}^{\infty}\frac{1}{[r]_{\q}![s]_{\q}!}\left(\sum_{l=0}^{s}\gauss{s}{l}_{\q}\q^{r(s-l)}a^{[r+l]}(z)b^{[s-l]}(z\q^{r+
l})\right)c(z\q^{s+r})z_0^{r}z_2^{s}\nonumber\\
&\qquad=\sum_{s= 0}^{\infty}\sum_{r=
0}^{\infty}\sum_{l=0}^{s}\frac{\q^{r(s-l)}}{[r]_{\q}![l]_{\q}![s-l]_{\q}!}a^{[r+l]}(z)b^{[s-l]}(z\q^{r+
l})c(z\q^{s+r})z_0^{r}z_2^{s};\label{last1}\\
&Y(a(z),z_0 +  z_2) Y(b(z),z_2) c(z) =Y(a(z),z_0
+
z_2)\sum_{t=0}^{\infty}\frac{1}{[t]_{\q}!}b^{[t]}(z)c(z\q^t)z_2^{t}\nonumber\\
&\qquad=\sum_{u=0}^{\infty}\sum_{t=0}^{\infty}\frac{1}{[u]_{\q}![t]_{\q}!}a^{[u]}(z)
b^{[t]}(z\q^u) c(z\q^{t+u})z_2^{t}(z_0 +z_2)^u\nonumber\\
&\qquad=\sum_{u=0}^{\infty}\sum_{t=0}^{\infty}\frac{1}{[u]_{\q}![t]_{\q}!}a^{[u]}(z)b^{[t]}(z\q^u) c(z\q^{t+u})z_2^{t}\left(\sum_{l=0}^{u}\gauss{u}{l}_{\q}z_0^{u-l}
z_2^l\right)
\nonumber\\
&\qquad=\sum_{u=0}^{\infty}\sum_{t=0}^{\infty}\sum_{l=0}^{u}\frac{\q^{t(u-l)}}{[u-l]_{\q}![l]_{\q}![t]_{\q}!}a^{[u]}(z)b^{[t]}(z\q^u) c(z\q^{t+u})z_0^{u-l}z_2^{t+l}.\label{last8}
\end{align}
Finally, by applying  substitutions
$r=u-l$ and  $s=t+l$ to \eqref{last8}
we get \eqref{last1}, thus proving \eqref{V4}.
\end{prf}

\begin{lem}\label{important1}
Every maximal quasi-commutative subspace of $\mathcal{E}(L)$ is closed.
\end{lem}

\begin{prf}
Let  $V$ be a maximal quasi-commutative subspace of $\mathcal{E}(L)$ and suppose
that there exist operators $a(z),b(z)\in V$ and an integer $r$ such $a(z)_r b(z)\notin V$. Then Lemma \ref{dong} implies that $V$ is not maximal. Contradiction.
\end{prf}

Notice that every maximal quasi-commutative subspace of $\mathcal{E}(L)$ contains identity $1_L$.
The following lemma is a  consequence of Theorem \ref{construction1} and
Lemma \ref{important1}.

\begin{lem}\label{important2}
Every maximal quasi-commutative subspace of $\mathcal{E}(L)$  is a nonlocal $\q$-vertex algebra.
\end{lem}

Finally, we formulate the main result of this section.

\begin{thm}\label{main}
Let $\mathcal{S}$ be a quasi-commutative subset of $\mathcal{E}(L)$. There exists a unique
smallest nonlocal $\q$-vertex algebra
$V\subset \mathcal{E}(L)$ such that $\mathcal{S}\subseteq V$.
\end{thm}

\begin{prf}
By Zorn's lemma, $\mathcal{S}$ is a subset of some maximal quasi-commutative subspace $W$
of $\mathcal{E}(L)$.  By
Lemma \ref{important2}, $W$ is a nonlocal $\q$-vertex algebra. Now, we  construct $V$ as the intersection of  all
 nonlocal $\q$-vertex algebras containing $\mathcal{S}$.
\end{prf}

Denote the resulted  nonlocal $\q$-vertex algebra $V$ by $\left<\mathcal{S}\right>$. We can use Theorem \ref{main} to construct modules for $\q$-vertex algebras as well.

\begin{kor}\label{maincor}
Let  $V\subset \mathcal{E}(L)$ be a nonlocal $\q$-vertex algebra and let $a(z)$ be an operator in $\om(K,L((z)))$, 
where $K$ is an arbitrary vector space over the field $\mathbb{F}$, such that 
$$b_1(z_1)\cdots b_m(z_m)a(z)\in \om(K,L((z_1,...,z_m,z)))$$
for all $b_1 (z),..., b_m(z)\in V$, $m\in\mathbb{Z}_{\geq 0}$.
There exists a unique smallest $V$-module $W$ such that $a(z)\in W$.
\end{kor}

Denote the resulted  $\left<\mathcal{S}\right>$-module by $\left<a(z)\right>$.

 \begin{rem}
By setting
$$a(z)\cdot b(z) = a(z)_{-1} b(z)\qquad \text{for all } a(z),b(z)\in  V,$$
the nonlocal $\q$-vertex algebra $V$, constructed in Theorem \ref{main}, becomes an associative algebra with unit $1$. 
Even though this algebra may not be commutative,
the underlying nonlocal $\q$-vertex algebraic structure can be recovered similarly as in Borcherds' construction (see \cite{Bor}).
More precisely, for $a(z),b(z)\in  V$ we have
$$ Y(a(z),z_0)b(z)=\sum_{r\geq 0}\frac{1}{[r]!}a^{[r]}(z) b(z\q^r)z_{0}^{r}=\left(\sum_{r\geq 0}\frac{1}{[r]!}a^{[r]}(z) z_{0}^{r}\right)b(z) =\left(e_{\q}^{z_0 \frac{\partial_{\q}}{\partial_{\q}z}}a(z)\right) \cdot b(z),$$
where, on the right side, $z_0$  denotes the right multiplication with variable $z_0$ satisfying
\eqref{asbefore}.
 \end{rem}

\section{Principal subspaces for \texorpdfstring{$U_{q}(\widehat{\mathfrak{sl}}_{2})$}{Uq(sl2)} and Frenkel-Jing operators}\label{FJsection} 

\subsection{Nonlocal \texorpdfstring{$\q$}{q}-vertex algebra \texorpdfstring{$\left<x(z)\right>$}{<x(z)>}} 

Set 
$\q =q^2 .$
In order to simplify our notation, we shall omit superscript ``$+$'' and write $x(z)$ instead of $x^{+}(z)$. This should not cause any confusion because the operator  $x^{-}(z)$ will not be used in the rest of the paper.

Recall  \eqref{remember}. By \eqref{r1} and Theorem \ref{main}  there exists a unique
smallest nonlocal $\q$-vertex algebra $\left<x(z)\right>\subset\mathcal{E}(L_0)$,
  with the vacuum vector $1=1_{L_0}:L_0 \to L_0$, which is
 generated by $\mathcal{S}=\left\{x(z)\right\}$.  
 
 \begin{rem}
 Although $a(z)_m b(z) =0$ for all $a(z),b(z)\in \left<x(z)\right>$, $m\geq 0$,  in general, the operators $Y(a(z),z_1)$ and $Y(b(z),z_2)$  do not commute.
 For example, consider the operators $x(z),x^{[1]}(z)=x(z)_{-2}1\in \left<x(z)\right>$. We have
\begin{align*}&\lim_{z_1,z_2 \to 0}\left[Y(x(z),z_1),Y(x^{[1]}(z),z_2)\right]1\\&\hspace{40pt}=x(z)_{-1}x^{[1]}(z)_{-1}1-x^{[1]}(z)_{-1} x(z)_{-1}1=\frac{x(z\q)x(z)}{z(1-\q)}\neq 0.\end{align*}
Specifically, the associative algebra $\left<x(z)\right>$ is not commutative.
 \end{rem}

Denote by $W_{0,\q}$ the subspace of $\left<x(z)\right>$ spanned by the operators
\begin{equation}\label{mono1}
 x(z)_{l_m}\ldots x(z)_{l_1}1 \in \left<x(z)\right>,\quad l_j \leq -1,\, j=1,2...,m,\, m\in\mathbb{Z}_{\geq 0}.
\end{equation}

\begin{rem}
Spaces $W_{0,\q}$ and $\left<x(z)\right>$ do not coincide. Set
 $$a(z)=(x(z)_{-2}1)_{-1}(x(z)_{-1}1)=\frac{x(z\q)-x(z)}{z(\q -1)}x(z)=x(z)_{-2}x(z\q^{-1})_{-1}1\in \left<x(z)\right>.$$
By using the techniques from the proof of Theorem \ref{maintheorem1} one can prove that $a(z)\notin W_{0,\q}$.
 \end{rem}
 
Equality  \eqref{r2} and Corollary \ref{maincor} imply that  there exists a unique smallest $\left<x(z)\right>$-module $\left<\mathcal{Y}(z)\right>\subset \om(L_0 ,L_1 ((z)))$ such that
$\mathcal{Y}(z)\in \left<\mathcal{Y}(z)\right>$. Denote by $W_{1,\q}$ the subspace of $\left<\mathcal{Y}(z)\right>$ spanned by the operators
\begin{equation}\label{mono2}
x(z)_{l_m}\ldots x(z)_{l_1}\mathcal{Y}(z) \in \left<\mathcal{Y}(z)\right>,\quad l_j \leq -1,\, j=1,2...,m,\, m\in\mathbb{Z}_{\geq 0}.
\end{equation}

\subsection{Basis for \texorpdfstring{$W_{i,\q}$}{Wiq}} \label{sec:100}
Spanning set \eqref{mono1} is an analogue of  spanning set \eqref{span1}. However, since the operator $x(z)$ is not commutative,
we can not assume that the sequence of indices $(l_m,...,l_1)$ in \eqref{mono1} is decreasing from right to left.
We proceed similarly as in Subsection \ref{sec:3}. First, we use nonlocal $\q$-vertex algebra structure to obtain some relations among the monomials in
\eqref{mono1}. Next, using these relations we  reduce the spanning set up to a basis.

The following Lemma is a special case of quantum integrability (cf. \cite{DM}). 

\begin{lem}\label{41}  
 On $\left<x(z)\right>$ we have
 \begin{align}
 & x(z)x(z) =x(z)_{-1} x(z)_{-1}  =  0;\label{rels1}\\
   & x(z)x^{[1]}(z) = x(z)_{-1} x(z)_{-2}  =  0;\label{rels2}\\
      & x^{[1]}(z)x(z\q) = x(z)_{-2} x(z)_{-1}  =  0.\label{rels2x}
 \end{align}
\end{lem}

\begin{prf}
The relations can be easily verified by applying $\displaystyle\lim_{\substack{z_1\to z\\z_2 \to z}}$ and  $\displaystyle\lim_{\substack{z_1\to z\\z_2 \to zq^2}}$ to \eqref{r1}.
\end{prf}

Applying  the $\q$-derivation $\frac{d^{n}_{\q}}{d_{\q} z^{n}}$ on \eqref{rels1},
$\q$-derivation
 $\frac{d^{n-1}_{\q}}{d_{\q} z^{n-1}} $ on \eqref{rels2}  and using 
Proposition \ref{qLeibniz} we get:

\begin{lem}\label{42}
 On $\left<x(z)\right>$ we have
 \begin{align}
 & \sum_{l=0}^{n}\gauss{n}{l}_{\q}x^{[l]}(z)x^{[n-l]}(z\q^l)=[n]_{\q}!\sum_{l=0}^{n} x(z)_{-l-1}x(z)_{-n+l-1} =0;\label{rels3}\\
   & \sum_{l=0}^{n-1}\gauss{n-1}{l}_{\q}x^{[l]}(z)x^{[n-l]}(z\q^l)=[n-1]_{\q}!\sum_{l=0}^{n-1}[n-l]_{\q} x(z)_{-l-1}x(z)_{-n+l-1} =0.\label{rels4}
 \end{align}
\end{lem}

\begin{rem}
Using the nonlocal $\q$-vertex algebra structure we were able to derive the  relations among operators $x(z)$,
 which are analogous to the relations among commutative operators $\widehat{x}(z)$ found earlier. More precisely,
 relations \eqref{rels1}--\eqref{rels2x} coincide with \eqref{niceexample1}, while \eqref{rels3} coincides with 
 \eqref{niceexample2}. Even though the analogue of \eqref{rels4} can be easily derived for the operators $\widehat{x}(z)$ as well,
 we did not needed it in  \Cref{komutativnisection} because, in that construction, we implicitly used another relation - commutativity.
  In the rest of this subsection, we use  similar techniques as in the proofs of Lemma \ref{second} and Theorem \ref{findthis}.
\end{rem}

 Set
 \begin{align}\label{W0qbaza}
\mathfrak{B}_{W_{0,\q}}=\bigg\{ x(z)_{l_m}\ldots x(z)_{l_1}1 \,\,\bigg| \,\,\,\bigg.
&l_1\leq -1\text{ and }\bigg.
 l_{r}\leq -3\\
\bigg. &\text{for all }l_r\in\mathbb{Z},\,r=2,3,\ldots,m,\,
m\in\mathbb{Z}_{\geq 0}\bigg\}.\nonumber
\end{align}

 \begin{lem}\label{spanningset}
 The set 
$\mathfrak{B}_{W_{0,\q}}$ spans  $W_{0,\q}$.
 \end{lem}
 
 \begin{prf}
  By  rewriting  \eqref{rels3} and \eqref{rels4} we get
 \begin{align}
  x(z)_{-1}x(z)_{-n-1}+x(z)_{-2}x(z)_{-n}&=-\sum_{l=2}^{n} x(z)_{-l-1}x(z)_{-n+l-1} ;\label{rel5}\\
    x(z)_{-1}x(z)_{-n-1}+\frac{[n-1]_{\q}}{[n]_{\q}} x(z)_{-2}x(z)_{-n}&=-\sum_{l=2}^{n-1}\frac{[n-l]_{\q}}{[n]_{\q}} x(z)_{-l-1}x(z)_{-n+l-1}.\label{rel6}
 \end{align}
 Consider the operator $a_j (z)=x(z)_{l_m}\ldots x(z)_{l_1}1$, where $l_j=-1,-2$ for some $j>1$ and $l_k\leq -3$ for $k\geq j+1$. 
 Relations \eqref{rel5} and \eqref{rel6} allow us to express $a_j(z)$ as a linear combination of the following (finite) families of operators:
 \begin{itemize}
  \item Operators $b(z)=x(z)_{k_m}\ldots x(z)_{k_1}1\in \mathfrak{B}_{W_{0,\q}}$;
  \item Operators $a_{j-1}(z)=x(z)_{t_m}\ldots x(z)_{t_1}1$, where $t_p\leq -3$ for $p\geq j$. 
 \end{itemize}
Applying the same argument on every operator  $a_{j-1}(z)$, then on every operator $a_{j-2}(z)$, etc.
we conclude that the original operator $a_{j}(z)$ can be expressed as a linear combination of the operators in $\mathfrak{B}_{W_{0,\q}}$
and the operators 
$x(z)_{u_m}\ldots x(z)_{u_1}1$ satisfying the following condition:
\begin{equation}\label{rel7}
 \textrm{there exists }l=0,1,...,m\textrm{ such that }u_j\leq -3\textrm{ for }j\geq l\textrm{ and }u_j=-1,-2\textrm{ for }j\leq l.
\end{equation}
Finally, relations \eqref{rels1}--\eqref{rels2x}, together with the fact that $x(z)_{-2}x(z)_{-2}1$ is proportional to $x(z)_{-3}x(z)_{-1}1$, imply that
in condition \eqref{rel7}
we can assume $l=0,1$. 
 \end{prf}
 
 Now, we proceed towards proving linear independence of the set $\mathfrak{B}_{W_{0,\q}}$.
 
 \begin{lem}\label{maksimalnimonomi}
  The set
  \begin{equation}\label{helpful}
   \left\{ :x(z\q^{l_m})\ldots x(z\q^{l_1}): \,\, \left| \,\,\right. l_m < l_{m-1} < \ldots < l_1,\, l_j\in\mathbb{Z},\,j=1,2,...,m,\, m\in\mathbb{Z}_{>0}\right\}
  \end{equation}
  is linearly independent.
 \end{lem}
 
 \begin{prf}
  Let
  $$\sum_{j=1}^{n} \alpha_j :x(z\q^{l_{m_j,j}})\ldots x(z\q^{l_{1,j}}): =0$$
  for some nonzero scalars $\alpha_j$, $j=1,2,...,n$. Without loss of generality we can assume that $n$ is a minimal integer for which such a nontrivial linear combination exists. 
  Since all the elements of the given set are nonzero, we have $n>1$.  Furthermore, we can assume that $m_j=m_k$ for all $j,k=1,2,...,n$. Indeed, if $m_j\neq m_k$ 
  for some $j,k=1,2,...,n$, by applying the operator $1\otimes \q^\alpha$ on the linear combination we get
    $$\sum_{j=1}^{n} \alpha_j \q^{2 m_j} :x(z\q^{l_{m_j,j}})\ldots x(z\q^{l_{1,j}}): =0,$$
    which, together with the original linear combination, contradicts to minimality of $n$. Therefore, by setting $m=m_j$, we have a linear combination
      $$\sum_{j=1}^{n} \alpha_j :x(z\q^{l_{m,j}})\ldots x(z\q^{l_{1,j}}): =0$$
      
   Naturally, we can also assume that all the summands $:x(z\q^{l_{m,j}})\ldots x(z\q^{l_{1,j}}):$ are different. 
   Then, there exist $p,r=1,2,...,n$ and $k=1,2,...,m$ such that
   $$l_{k,p}\neq l_{j,r}\quad \textrm{for }j=1,2,...,m.$$ 
 Recall \eqref{komutiranjesE-}. Multiplying the linear combination from the right side with $\mathcal{E}^{-}(z_1\q^{l_{k,p}})$ and from the left side with $\mathcal{E}^{-}(z_1\q^{l_{k,p}})^{-1}$ and 
 applying the limit 
 $\lim_{z_1\to z}$ we get
 $$\sum_{j=1}^{n} \beta_j :x(z\q^{l_{m,j}})\ldots x(z\q^{l_{1,j}}): =0$$
 for some scalars $\beta_j$ such that
 $\beta_p = 0$ and $\beta_r \neq 0$.
This is a contradiction to minimality of $n$.
 \end{prf}

 \begin{thm}\label{maintheorem1}
  The set $\mathfrak{B}_{W_{0,\q}}$ forms a basis for ${W_{0,\q}}$.
 \end{thm}

 \begin{prf}
 Notice that every element of $\mathfrak{B}_{W_{0,\q}}$ is nonzero because it can be written as a $\mathbb{C}(q^{1/2})[z^{\pm 1}]$-linear 
  combination of linearly independent elements \eqref{helpful}. 
  For example, the element
  $x(z)_{-3} x(z)_{-2} 1$ is a 
  $\mathbb{C}(q^{1/2})[z^{\pm 1}]$-linear combination of
  $$:x(z)x(z\q^2):,\quad :x(z)x(z\q^3):,\quad :x(z\q)x(z\q^3):.$$
  
  Let
  $$\sum_{j=1}^{n} \alpha_j x(z)_{l_{m_j ,j}}\ldots x(z)_{l_{1 ,j}}1 =0 $$
  be a linear combination such that all the summands $x(z)_{l_{m_j ,j}}\ldots x(z)_{l_{1 ,j}}1$
  are different elements of $\mathfrak{B}_{W_{0,\q}}$ and  $\alpha_j \neq 0$ for $j=1,2,...,n$.
Without loss of generality we can assume that $n$ is a minimal integer for which such a nontrivial linear combination exists.  
Using the same argument as in the proof of Lemma \ref{maksimalnimonomi} we can assume that $m_j=m_k$ for all $j,k=1,2,...,n$.
Therefore, by setting $m=m_j$, we have a linear combination
 $$\sum_{j=1}^{n} \alpha_j x(z)_{l_{m ,j}}\ldots x(z)_{l_{1 ,j}}1 =0.$$
 
 For $k=1,2,...,m$ set
 \begin{align}
  &D_k (j)=D_k (x(z)_{l_{m ,j}}\ldots x(z)_{l_{1 ,j}}1)=\sum_{p=1}^{k} l_{p,j}+k-1;\label{lj}\\
  &D_{k}^{*}(j)=D_k (j) -D_m (j).
 \end{align}
 Recall that every $x(z)_{l_{m ,j}}\ldots x(z)_{l_{1 ,j}}1$ can be written as a $\mathbb{C}(q^{1/2})[z^{\pm 1}]$-linear combination of elements 
 \eqref{helpful}.
 The integer $-D_m(j)-1$ is the maximal integer for which at least one element of this linear combination contains $x(z\q^{-D_m(j)-1})$.
 For example,
$$D_2 (x(z)_{-3} x(z)_{-2} 1)=(-3)+(-2)+2-1=-4.$$
 Also,  applying the operator $x(z)_{l_{m ,j}}\ldots x(z)_{l_{1 ,j}}1$ on the vector $1\otimes 1 \in L_0$ we get a Laurent  series
 $a_j (z) \in z^{D_m(j)+1+m(m-1)} L_0 [[z]]$. The lowest power of the variable $z$ in $a_j (z)$ is exactly $z^{D_m(j)+1+m(m-1)}$.
 
 We can express every summand $x(z)_{l_{m ,j}}\ldots x(z)_{l_{1 ,j}}1$ in the following way:
 \begin{align*}
  x(z)_{l_{m ,j}}\ldots x(z)_{l_{1 ,j}}1\, =\, \beta_j b_j (z) z^{m(m-1)}\, +\,\textrm{some other summands},
 \end{align*}
 where
 $$b_j (z)=\frac{:x(z\q^{D_{m}^{*}(j)})x(z\q^{D_{m-1}^{*}(j)})...x(z\q^{D_{1}^{*}(j)}):}{z^{-D_m(j)-1}}$$
 and $\beta_j$ is a nonzero scalar.
  Lemma \ref{maksimalnimonomi} implies that the set $\left\{b_1 (z),...,b_n (z)\right\}$ is linearly independent, so we conclude that the set 
  $\left\{ x(z)_{l_{m ,j}}\ldots x(z)_{l_{1 ,j}}1\, \left| \right.\, j=1,2,...,n\right\}$ is linearly independent as well. Hence, $\alpha_j =0$ for
  $j=1,2,..,n$, so  Lemma \ref{spanningset} implies the statement of the theorem.
 \end{prf}

 In the rest of this subsection, we briefly explain how to construct a similar basis for ${W_{1,\q}}$. 
 First, we notice that Lemmas \ref{41} and \ref{42} can be applied in this case as well.  Next, we define
  \begin{align}\label{W1qbaza}
\mathfrak{B}_{W_{1,\q}}=\bigg\{ x(z)_{l_m}\ldots x(z)_{l_1}\mathcal{Y}(z) \,\,\bigg| \,\,\,\bigg.
&l_1\leq -2\text{ and }\bigg.
 l_{r}\leq -3\\
\bigg. &\text{for all }l_r\in\mathbb{Z},\,r=2,3,\ldots,m,\,
m\in\mathbb{Z}_{\geq 0}\bigg\}.\nonumber
\end{align}
Condition $l_1\leq -2$ is a consequence of  \eqref{r2} because by applying the limit $\lim_{z_2 ,z_1 \to z}$ to this relation
 we get
 $$x(z)_{-1}\mathcal{Y}(z)=x(z)\mathcal{Y}(z)=0.$$
Therefore, 
 
 \begin{lem}\label{obviously}
 The set $\mathfrak{B}_{W_{1,\q}}$ spans  $W_{1,\q}$.
 \end{lem}
 
 The following result can be proved in the same way as Lemma \ref{maksimalnimonomi}.
  \begin{lem}\label{maksimalnimonomi2}
  For every integer $l_0$ the set
  \begin{equation*}
   S=\left\{ :x(z\q^{l_m})\ldots x(z\q^{l_1})\mathcal{Y}(z\q^{l_0}): \,\, \left| \,\,\right. l_m < \ldots < l_1 < l_0,\, l_j\in\mathbb{Z},\,j=1,2,...,m,\, m\in\mathbb{Z}_{>0}\right\}
  \end{equation*}
  is linearly independent.
 \end{lem}
 
 By right
 multiplying  $:x(z\q^{l_m})\ldots x(z\q^{l_1})\mathcal{Y}(z\q^{l_0}):\in S$   by $E_{+}^{+}(-a ,z_1 \q^{l})^{-1}$, 
 left multiplying by $E_{+}^{+}(-a ,z_1 \q^{l})$ and 
 applying the limit 
 $\lim_{z_1\to z}$ we get (recall \eqref{komutiranjesE2})
 $$\beta \delta _{ll_0}:x(z\q^{l_m})\ldots x(z\q^{l_1})\mathcal{Y}(z\q^{l_0}):$$
 for some nonzero scalar $\beta$. Hence, using the same technique as in the proof of Lemma \ref{maksimalnimonomi} we can prove:
 
  \begin{lem}\label{maksimalnimonomi3}
  The set
  \begin{equation}
   \left\{ :x(z\q^{l_m})\ldots x(z\q^{l_1})\mathcal{Y}(z\q^{l_0}): \,\, \left| \,\,\right. l_m < \ldots < l_1 < l_0,\, l_j\in\mathbb{Z},\,j=0,1,...,m,\, m\in\mathbb{Z}_{>0}\right\}
  \end{equation}
  is linearly independent.
 \end{lem}
 
  \begin{thm}\label{maintheorem2}
  The set $\mathfrak{B}_{W_{1,\q}}$ forms a basis for ${W_{1,\q}}$.
 \end{thm}
 
 \begin{prf}
 The proof goes analogously to the proof of Theorem \ref{maintheorem1}, but, at the end, we  express 
 every summand $x(z)_{l_{m ,j}}\ldots x(z)_{l_{1 ,j}}\mathcal{Y}(z)$ in the following way:
 \begin{align*}
  x(z)_{l_{m ,j}}\ldots x(z)_{l_{1 ,j}}\mathcal{Y}(z)\, =\, \beta_j b_j (z) z^{m^2}\, +\,\textrm{some other summands},
 \end{align*}
 where
 $$b_j (z)=\frac{:x(z\q^{D_{m}^{*}(j)})x(z\q^{D_{m-1}^{*}(j)})...x(z\q^{D_{1}^{*}(j)})\mathcal{Y}(z\q^{-D_m(j)-1}):}{z^{-D_m(j)-1}} $$
 and $\beta_j$ is a nonzero scalar. Now,  Lemma \ref{maksimalnimonomi3} implies linear independence, so the statement of the theorem follows from Lemma \ref{obviously}.
 \end{prf}

\subsection{On the sum side of Rogers-Ramanujan identities.} 
 
The bases $\widehat{\mathfrak{B}}_{W(\Lambda_i)}$, $i=0,1$, found in  \Cref{sec:2}, as well as the bases $\widehat{\mathfrak{B}}_{W_i}$, $i=0,1$,
found in \Cref{sec:3}, correspond to the sum side of Rogers-Ramanujan identities
\begin{align}\label{rrrrrrrrrr}
\prod_{r \geq 0}\frac{1}{(1-q^{5r+1+i})(1-q^{5r+4-i})}&= \sum_{r\geq 0}\frac{q^{r^2+ir}}{(1-q)(1-q^{2})\cdots (1-q^r)}.
\end{align}
 However, the bases $\mathfrak{B}_{W_{i,\q}}$, $i=0,1$, found in  \Cref{sec:100} are defined in terms of different  conditions among monomial indices, 
so it is not clear whether they are  related to \eqref{rrrrrrrrrr}. In this subsection, we show that they in fact correspond to the same identities.

Recall  \eqref{lj} and for an arbitrary basis element $$a(z)=x(z)_{l_{m }}\ldots x(z)_{l_{1 }}1\in \mathfrak{B}_{W_{0,\q}}\quad\textrm{or}\quad 
a(z)=x(z)_{l_{m }}\ldots x(z)_{l_{1 }}\mathcal{Y}(z)\in \mathfrak{B}_{W_{1,\q}}$$ 
define
\begin{equation*}
  \deg_{\q} a(z) =-\sum_{k=1}^{m} D_k (x(z)_{l_{m }}\ldots x(z)_{l_{1 }}1).
\end{equation*}
For an integer $r\geq 0$ and $i=0,1$ set
$$\left(W_{i,\q}\right)_r =\left\{a(z)\in W_{i,\q}\,\left|\right. \, \textstyle\deg_{\q} a(z)=r\right\}.$$
We have the direct sum decomposition
$$W_{i,\q} = \bigoplus_{r\geq 0} \left(W_{i,\q}\right)_r$$
and the subspaces $\left(W_{i,\q}\right)_r$ are finite-dimensional. Hence, we can introduce a character $\ch_{\q}$ by
$$\textstyle\ch_{\q}W_{i,\q}=\displaystyle\sum_{r\geq 0} \dim \left(W_{i,\q}\right)_r q^r .$$
Naturally, parameters $\q$ and $q$ in the above formula  are not related.
\begin{thm}\label{rr}
For $i=0,1$ we have
 $$\textstyle\ch_{\q}\displaystyle W_{i,\q}=\sum_{r\geq 0}\frac{q^{r^2+ir}}{(1-q)(1-q^{2})\cdots (1-q^r)}.$$
\end{thm}

\begin{figure}[H]
\begin{tikzpicture}[scale=1.1]
\tikzstyle{every node}=[font=\scriptsize]
\draw  (0,1) -- (0,0) -- (1,0) -- (1,1) -- (0,1);
\node at (0.5,0.5) {$x(z)$};

\draw  (3,1) -- (3,0) -- (4,0) -- (4,1) -- (3,1);
\node at (3.5,0.5) {$x(z)$};
\draw  (3,2) -- (3,1) -- (4,1) -- (4,2) -- (3,2);
\node at (3.5,1.5) {$x(z\q)$};

\draw  (6,1) -- (6,0) -- (7,0) -- (7,1) -- (6,1);
\node at (6.5,0.5) {$x(z)$};
\draw  (6,2) -- (6,1) -- (7,1) -- (7,2) -- (6,2);
\node at (6.5,1.5) {$x(z\q)$};
\draw  (6,3) -- (6,2) -- (7,2) -- (7,3) -- (6,3);
\node at (6.5,2.5) {$x(z\q^2)$};

\node at (0.5,-0.5) {$x(z)_{-1}1$};
\node at (3.5,-0.5) {$x(z)_{-2}1$};
\node at (6.5,-0.5) {$x(z)_{-3}1$};
\end{tikzpicture}\caption{Elements of basis $\mathfrak{B}_{W_{0,\q}}$ represented by Young diagrams}
\label{pic1}
\end{figure}
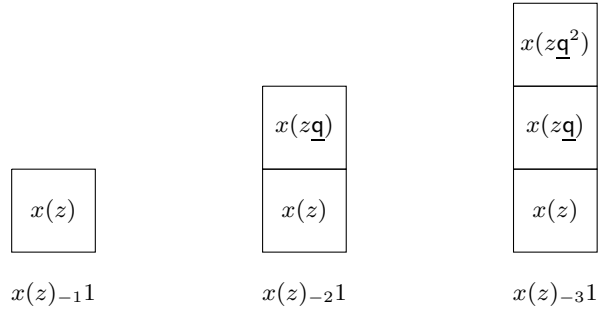

\begin{figure}[H]
\begin{tikzpicture}[scale=1.1]
\tikzstyle{every node}=[font=\scriptsize]
\draw  (1,3) -- (1,2) -- (2,2) -- (2,3) -- (1,3);
\node at (1.5,2.5) {$x(z\q^2)$};
\draw  (0,1) -- (0,0) -- (1,0) -- (1,1) -- (0,1);
\node at (0.5,0.5) {$x(z)$};
\draw  (0,2) -- (0,1) -- (1,1) -- (1,2) -- (0,2);
\node at (0.5,1.5) {$x(z\q)$};
\draw  (0,3) -- (0,2) -- (1,2) -- (1,3) -- (0,3);
\node at (0.5,2.5) {$x(z\q^2)$};

\draw  (4,1) -- (4,0) -- (5,0) -- (5,1) -- (4,1);
\node at (4.5,0.5) {$x(z)$};
\draw  (4,2) -- (4,1) -- (5,1) -- (5,2) -- (4,2);
\node at (4.5,1.5) {$x(z\q)$};
\draw  (4,3) -- (4,2) -- (5,2) -- (5,3) -- (4,3);
\node at (4.5,2.5) {$x(z\q^2)$};
\draw  (5,3) -- (5,2) -- (6,2) -- (6,3) -- (5,3);
\node at (5.5,2.5) {$x(z\q^2)$};
\draw  (5,4) -- (5,3) -- (6,3) -- (6,4) -- (5,4);
\node at (5.5,3.5) {$x(z\q^3)$};

\draw  (8,1) -- (8,0) -- (9,0) -- (9,1) -- (8,1);
\node at (8.5,0.5) {$x(z)$};
\draw  (8,2) -- (8,1) -- (9,1) -- (9,2) -- (8,2);
\node at (8.5,1.5) {$x(z\q)$};
\draw  (8,3) -- (8,2) -- (9,2) -- (9,3) -- (8,3);
\node at (8.5,2.5) {$x(z\q^2)$};
\draw  (9,3) -- (9,2) -- (10,2) -- (10,3) -- (9,3);
\node at (9.5,2.5) {$x(z\q^2)$};
\draw  (9,4) -- (9,3) -- (10,3) -- (10,4) -- (9,4);
\node at (9.5,3.5) {$x(z\q^3)$};
\draw  (9,5) -- (9,4) -- (10,4) -- (10,5) -- (9,5);
\node at (9.5,4.5) {$x(z\q^4)$};

\node at (1,-0.5) {$x(z)_{-3}x(z)_{-1}1$};
\node at (5,-0.5) {$x(z)_{-3}x(z)_{-2}1$};
\node at (9,-0.5) {$x(z)_{-3}x(z)_{-3}1$};

\end{tikzpicture}\caption{Elements $x(z)_{-3}x(z)_{-1}1, x(z)_{-3}x(z)_{-2}1, x(z)_{-3}x(z)_{-3}1\in \mathfrak{B}_{W_{0,\q}}$}
\label{pic2}
\end{figure}
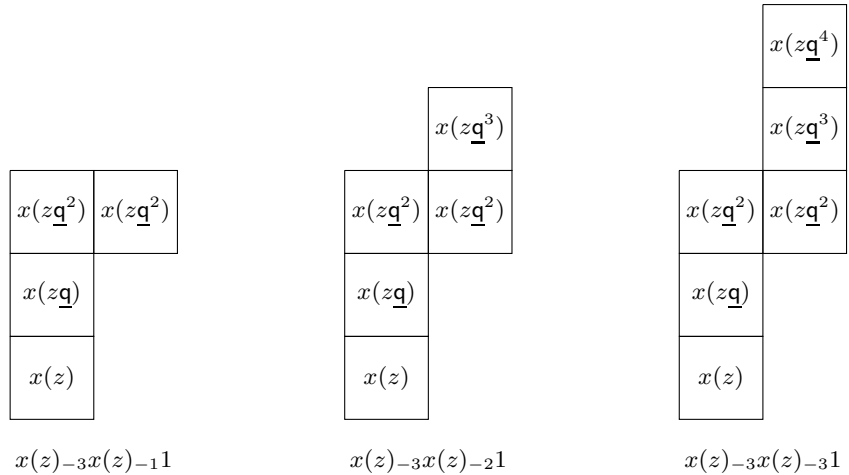

The statement of the Theorem can be easily verified by a convenient visualization of the basis elements.
We can represent elements of $\mathfrak{B}_{W_{i,\q}}$ by diagrams. For example, elements 
$$x(z)_{-1}1=x(z),\quad x(z)_{-2}1=\frac{x(z\q) -x(z)}{z(\q -1)},\quad  x(z)_{-3}1=\frac{x(z\q^2) -(1+\q)x(z\q)+\q x(z)}{\q z^2(\q -1)^2}$$
are shown in Figure \ref{pic1}.
We can join multiple diagrams, thus obtaining more complicated elements. For example,
elements
\begin{align*}
&x(z)_{-3}x(z)_{-1}1=\frac{x(z\q^2) -(1+\q)x(z\q)+\q x(z)}{\q z^2(\q -1)^2}\cdot x(z\q^2);\\
&x(z)_{-3}x(z)_{-2}1= \frac{x(z\q^2) -(1+\q)x(z\q)+\q x(z)}{\q z^2(\q -1)^2}\cdot \frac{x(z\q^3) -x(z\q^2)}{z\q^2 (\q -1)};\\
&x(z)_{-3}x(z)_{-3}1=\frac{x(z\q^2) -(1+\q)x(z\q)+\q x(z)}{\q z^2(\q -1)^2}\cdot \frac{x(z\q^4) -(1+\q)x(z\q^3)+\q x(z\q^2)}{\q^5 z^2(\q -1)^2}
\end{align*}
are represented in Figure \ref{pic2}.

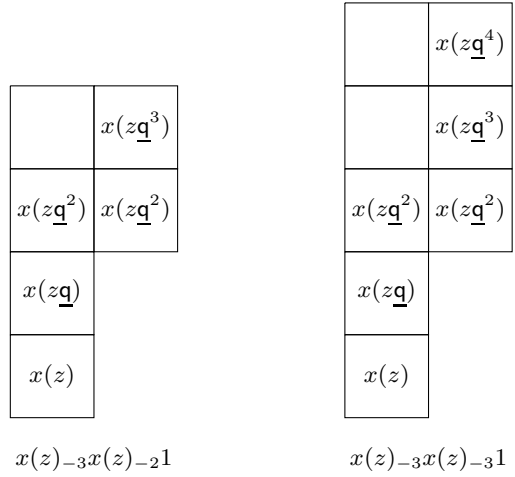
\begin{figure}[H]
\begin{tikzpicture}[scale=1.1]
\tikzstyle{every node}=[font=\scriptsize]
\draw  (4,1) -- (4,0) -- (5,0) -- (5,1) -- (4,1);
\node at (4.5,0.5) {$x(z)$};
\draw  (4,2) -- (4,1) -- (5,1) -- (5,2) -- (4,2);
\node at (4.5,1.5) {$x(z\q)$};
\draw  (4,3) -- (4,2) -- (5,2) -- (5,3) -- (4,3);
\node at (4.5,2.5) {$x(z\q^2)$};
\draw  (4,4) -- (4,3) -- (5,3) -- (5,4) -- (4,4);
\draw  (5,3) -- (5,2) -- (6,2) -- (6,3) -- (5,3);
\node at (5.5,2.5) {$x(z\q^2)$};
\draw  (5,4) -- (5,3) -- (6,3) -- (6,4) -- (5,4);
\node at (5.5,3.5) {$x(z\q^3)$};

\draw  (8,1) -- (8,0) -- (9,0) -- (9,1) -- (8,1);
\node at (8.5,0.5) {$x(z)$};
\draw  (8,2) -- (8,1) -- (9,1) -- (9,2) -- (8,2);
\node at (8.5,1.5) {$x(z\q)$};
\draw  (8,3) -- (8,2) -- (9,2) -- (9,3) -- (8,3);
\node at (8.5,2.5) {$x(z\q^2)$};
\draw  (8,4) -- (8,3) -- (9,3) -- (9,4) -- (8,4);
\draw  (8,5) -- (8,4) -- (9,4) -- (9,5) -- (8,5);
\draw  (9,3) -- (9,2) -- (10,2) -- (10,3) -- (9,3);
\node at (9.5,2.5) {$x(z\q^2)$};
\draw  (9,4) -- (9,3) -- (10,3) -- (10,4) -- (9,4);
\node at (9.5,3.5) {$x(z\q^3)$};
\draw  (9,5) -- (9,4) -- (10,4) -- (10,5) -- (9,5);
\node at (9.5,4.5) {$x(z\q^4)$};

\node at (5,-0.5) {$x(z)_{-3}x(z)_{-2}1$};
\node at (9,-0.5) {$x(z)_{-3}x(z)_{-3}1$};

\end{tikzpicture}\caption{Completed diagrams for $x(z)_{-3}x(z)_{-2}1, x(z)_{-3}x(z)_{-3}1\in \mathfrak{B}_{W_{0,\q}}$}
\label{pic3}
\end{figure}

The diagrams corresponding to the elements of $\mathfrak{B}_{W_{i,\q}}$ can be ''completed" by adding the minimal (nonnegative) number of
''empty boxes", thus getting a Young diagram. In Figure \ref{pic3} we can see  completed diagrams for
$x(z)_{-3}x(z)_{-2}1,x(z)_{-3}x(z)_{-3}1\in \mathfrak{B}_{W_{0,\q}}$. Notice that the diagram for $x(z)_{-3}x(z)_{-1}1$ is already a Young diagram, 
so we do not need to add any empty boxes.

Such a completed diagram is equal to the Young diagram of some  vector  $b\in\widehat{\mathfrak{B}}_{W(\Lambda_i)}$. More precisely, the completed diagram of the element
$$a_0 (z)=x(z)_{l_{m }}\ldots x(z)_{l_{1 }}1\in \mathfrak{B}_{W_{0,\q}}\quad\textrm{or}\quad 
a_1 (z)=x(z)_{l_{m }}\ldots x(z)_{l_{1 }}\mathcal{Y}(z)\in \mathfrak{B}_{W_{1,\q}}$$ is equal to the Young diagram
of the element
$$b_0 =\widehat{x}(D_m)\cdots\widehat{x}(D_1)v_{\Lambda_0}\in\widehat{\mathfrak{B}}_{W(\Lambda_0)} \quad\textrm{or}\quad b_1 = \widehat{x}(D_m)\cdots\widehat{x}(D_1)v_{\Lambda_1}\in\widehat{\mathfrak{B}}_{W(\Lambda_1)}$$
respectively, where $D_j=D_j(a_0 (z))$  for $j=1,2,...,m$.
Furthermore, we have
$$\deg_{\q} a_i (z) = -(D_m +...+D_2+D_1) = \deg b_i,\qquad i=0,1.$$

Denote by $\mathcal{D}(\mathfrak{B})$ a family of diagrams corresponding to the elements of a basis $\mathfrak{B}$.
By considering defining conditions \eqref{W0qbaza} and \eqref{W1qbaza} for $\mathfrak{B}_{W_{i,\q}}$   and 
difference conditions \eqref{baza} for $\widehat{\mathfrak{B}}_{W(\Lambda_i)}$,
we see that the above ''completion of diagrams" establishes a bijection 
$\mathcal{D}(\mathfrak{B}_{W_{i,\q}})\to \mathcal{D}(\widehat{\mathfrak{B}}_{W(\Lambda_i)})$ for $i=0,1$,
so the statement of Theorem \ref{rr} clearly follows.

\begin{figure}
\begin{tikzpicture}[scale=1.0]
\tikzstyle{every node}=[font=\tiny]
\draw  (0,1) -- (0,0) -- (1,0) -- (1,1) -- (0,1);
\node at (0.5,0.5) {$x(z)$};
\draw  (0,2) -- (0,1) -- (1,1) -- (1,2) -- (0,2);
\node at (0.5,1.5) {$x(z\q)$};
\draw  (0,3) -- (0,2) -- (1,2) -- (1,3) -- (0,3);
\node at (0.5,2.5) {$x(z\q^2)$};
\draw  (0,4) -- (0,3) -- (1,3) -- (1,4) -- (0,4);
\node at (0.5,3.5) {$x(z\q^3)$};

\draw  (1,4) -- (1,3) -- (2,3) -- (2,4) -- (1,4);
\node at (1.5,3.5) {$x(z\q^3)$};
\draw  (1,5) -- (1,4) -- (2,4) -- (2,5) -- (1,5);
\node at (1.5,4.5) {$x(z\q^4)$};
\draw  (1,6) -- (1,5) -- (2,5) -- (2,6) -- (1,6);
\node at (1.5,5.5) {$x(z\q^5)$};

\draw  (2,6) -- (2,5) -- (3,5) -- (3,6) -- (2,6);
\node at (2.5,5.5) {$x(z\q^5)$};
\draw  (2,7) -- (2,6) -- (3,6) -- (3,7) -- (2,7);
\node at (2.5,6.5) {$x(z\q^6)$};
\draw  (2,8) -- (2,7) -- (3,7) -- (3,8) -- (2,8);
\node at (2.5,7.5) {$x(z\q^7)$};
\draw  (2,9) -- (2,8) -- (3,8) -- (3,9) -- (2,9);
\node at (2.5,8.5) {$x(z\q^8)$};
\draw  (2,10) -- (2,9) -- (3,9) -- (3,10) -- (2,10);
\node at (2.5,9.5) {$x(z\q^9)$};

\draw  (3,10) -- (3,9) -- (4,9) -- (4,10) -- (3,10);
\node at (3.5,9.5) {$x(z\q^9)$};
\draw  (3,11) -- (3,10) -- (4,10) -- (4,11) -- (3,11);
\node at (3.5,10.5) {$x(z\q^{10})$};

\draw  (8,1) -- (8,0) -- (9,0) -- (9,1) -- (8,1);
\draw  (8,2) -- (8,1) -- (9,1) -- (9,2) -- (8,2);
\draw  (8,3) -- (8,2) -- (9,2) -- (9,3) -- (8,3);
\draw  (8,4) -- (8,3) -- (9,3) -- (9,4) -- (8,4);
\draw  (8,5) -- (8,4) -- (9,4) -- (9,5) -- (8,5);
\draw  (8,6) -- (8,5) -- (9,5) -- (9,6) -- (8,6);
\draw  (8,7) -- (8,6) -- (9,6) -- (9,7) -- (8,7);
\draw  (8,8) -- (8,7) -- (9,7) -- (9,8) -- (8,8);
\draw  (8,9) -- (8,8) -- (9,8) -- (9,9) -- (8,9);
\draw  (8,10) -- (8,9) -- (9,9) -- (9,10) -- (8,10);
\draw  (8,11) -- (8,10) -- (9,10) -- (9,11) -- (8,11);

\draw  (9,4) -- (9,3) -- (10,3) -- (10,4) -- (9,4);
\draw  (9,5) -- (9,4) -- (10,4) -- (10,5) -- (9,5);
\draw  (9,6) -- (9,5) -- (10,5) -- (10,6) -- (9,6);
\draw  (9,7) -- (9,6) -- (10,6) -- (10,7) -- (9,7);
\draw  (9,8) -- (9,7) -- (10,7) -- (10,8) -- (9,8);
\draw  (9,9) -- (9,8) -- (10,8) -- (10,9) -- (9,9);
\draw  (9,10) -- (9,9) -- (10,9) -- (10,10) -- (9,10);
\draw  (9,11) -- (9,10) -- (10,10) -- (10,11) -- (9,11);

\draw  (10,6) -- (10,5) -- (11,5) -- (11,6) -- (10,6);
\draw  (10,7) -- (10,6) -- (11,6) -- (11,7) -- (10,7);
\draw  (10,8) -- (10,7) -- (11,7) -- (11,8) -- (10,8);
\draw  (10,9) -- (10,8) -- (11,8) -- (11,9) -- (10,9);
\draw  (10,10) -- (10,9) -- (11,9) -- (11,10) -- (10,10);
\draw  (10,11) -- (10,10) -- (11,10) -- (11,11) -- (10,11);

\draw  (11,10) -- (11,9) -- (12,9) -- (12,10) -- (11,10);
\draw  (11,11) -- (11,10) -- (12,10) -- (12,11) -- (11,11);

\node[font=\normalsize] at (2,-1) {$x(z)_{-4}x(z)_{-3}x(z)_{-5}x(z)_{-2}\mathcal{Y}(z)$};
\node[font=\normalsize] at (10,-1) {$\widehat{x}(-11)\widehat{x}(-8)\widehat{x}(-6)\widehat{x}(-2)v_{\Lambda_1}$};
\node[font=\normalsize] at (2,12) {$\mathcal{D}(\mathfrak{B}_{W_{1,\q}})$};
\node[font=\normalsize] at (10,12) {$\mathcal{D}(\widehat{\mathfrak{B}}_{W(\Lambda_1)})$};
\draw[->]  (5,12) -- (7,12);
\draw[|->] (5,5.5) -- (7,5.5);  

\end{tikzpicture}\caption{Bijection $\mathcal{D}(\mathfrak{B}_{W_{1,\q}})\to \mathcal{D}(\widehat{\mathfrak{B}}_{W(\Lambda_1)})$}
\label{pic4}
\end{figure}
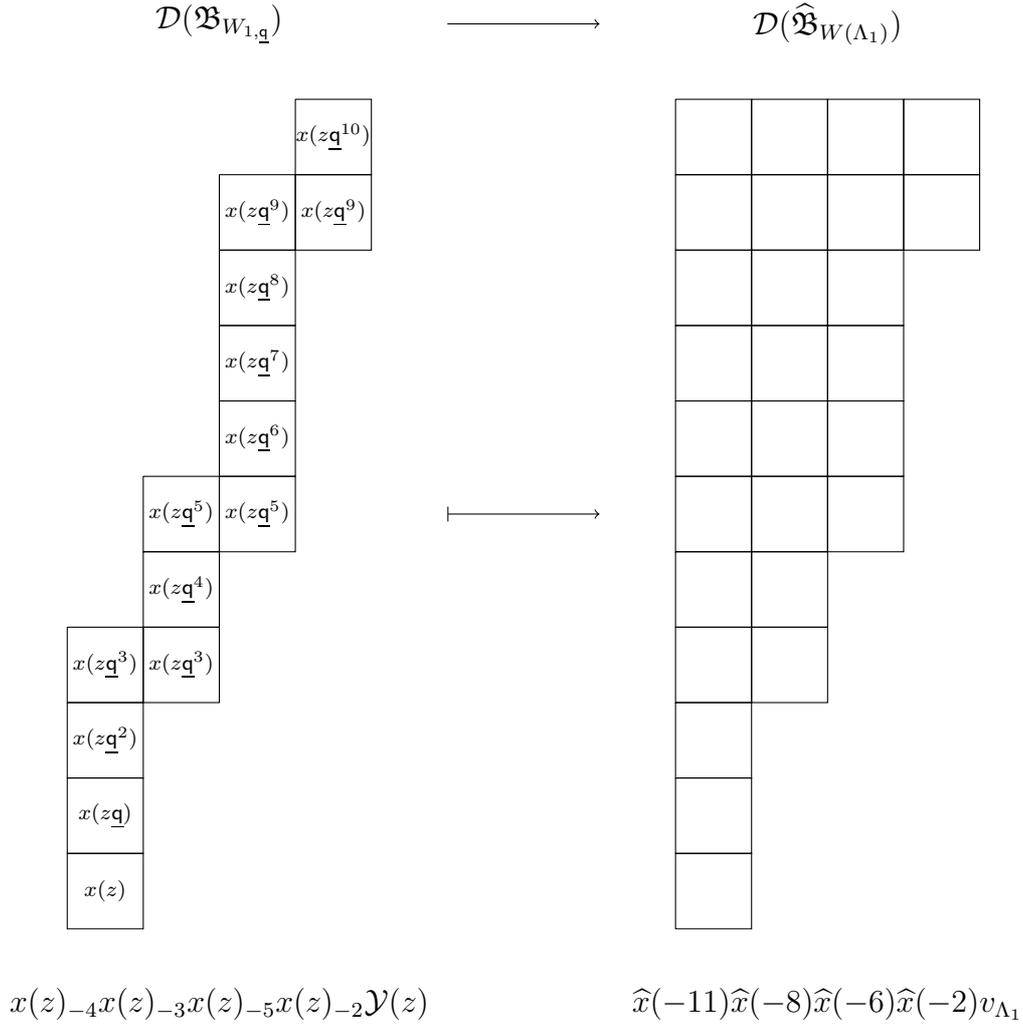

\begin{ex}
The completed diagram of $x(z)_{-4}x(z)_{-3}x(z)_{-5}x(z)_{-2}\mathcal{Y}(z)\mathfrak{B}_{W_{1,\q}}$ is equal to the Young diagram
of $\widehat{x}(-11)\widehat{x}(-8)\widehat{x}(-6)\widehat{x}(-2)v_{\Lambda_1}\in\widehat{\mathfrak{B}}_{W(\Lambda_1)}$ (Figure \ref{pic4}).
\end{ex}

\section{Relations among quantum quasi-particles}\label{zadnjisection}

\subsection{Quasi-particles} 
Let $V$ be an integrable highest weight $\widehat{\mathfrak{sl}}_{n+1}$-module of level $c$.
In \cite{G}, Georgiev defined quasi-particles of charge $m\in\mathbb{Z}_{>0}$ as coefficients
of the operator
\begin{equation}\label{intro1}
x_{m\alpha}(z)=\sum_{r\in\mathbb{Z}} x_{m\alpha}(r)z^{-r-m}=\left(x_{\alpha}(z)\right)^m \in\om(V,V((z))),
\end{equation}
where $\alpha$ is a simple root of the simple  Lie algebra $\mathfrak{sl}_{n+1}$.
Furthermore, for positive integers $m$ and $k$, $m\leq k$, he found $2m$ relations among
$x_{m\alpha}(z)$ and $x_{k\alpha}(z)$. Using the relations he expressed, for an arbitrary $v\in V$,  $2m$ vectors
\begin{equation*}
x_{m\alpha}(r+l)x_{k\alpha}(N-(r+l))v,\qquad
l=0,1,\ldots,2m-1,
\end{equation*}
as a linear combination of  vectors
\begin{equation*}
x_{m\alpha}(s)x_{k\alpha}(N-s)v,\qquad
s\in\mathbb{Z}\setminus\left\{r,r+1,\ldots,r+2m-1\right\},
\end{equation*}
 and vectors
\begin{equation*}
x_{m'\alpha}(t)x_{k'\alpha}(N-t)v,\qquad
0\leq m'<m,\quad m+k=m'+k',\quad t\in\mathbb{Z}.
\end{equation*}
These relations can be written as
\begin{equation}\label{intro2}
\left(\frac{d^l}{dz^l}x_{m\alpha}(z)\right)x_{k\alpha}(z)=\left(A_l
(z) +B_l (z) \frac{d^l}{dz^l}\right)x_{(k+1)\alpha}(z),
\end{equation}
 where $l=0,1,\ldots,2m-1$ and $A_l (z)$, $B_l (z)$ are some formal power series
(see \cite[Lemma 4.2.]{JP} for details).
The above formulae have the following vertex-operator theoretic interpretation. 
Relation \eqref{intro1} can be written as
\begin{equation}\label{intro3}
x_{m\alpha}(z)=\underbrace{x_{\alpha}(z)_{-1}\ldots
x_{\alpha}(z)_{-1} x_{\alpha}(z)_{-1}}_{m}1,
\end{equation}
while the left side in \eqref{intro2} equals
\begin{equation}\label{intro4}
x_{m\alpha}(z)_{-l-1} x_{k\alpha}(z),\qquad l=0,1,\ldots,2m-1
\end{equation}
(up to a multiplication scalar). In the next subsection we provide a similar
vertex-operator theoretic interpretation of the relations among quantum quasi-particles.

\subsection{Quantum quasi-particles} 
As in the previous section, set $\q =q^2$.
Let $L$ be an arbitrary  integrable highest weight $U_q(\widehat{\mathfrak{sl}}_{n+1})$-module of level $c\in\mathbb{Z}_{> 0}$.
We denote by $\bar{x}^{+}_{\alpha_{i}}(z)\in\mathcal{E}(L)$,  $i=1,2,...,n$, Ding-Feigin operators introduced in \cite{DF}.
They satisfy
\begin{align}
&[\bar{x}^{+}_{\alpha_{i}}(z_1),\bar{x}^{+}_{\alpha_{i}}(z_2)]=0,\label{dfkomutativnost}\\
&\bar{x}^{+}_{\alpha_{i}}(z) \bar{x}^{+}_{\alpha_{i}}(zq^2)\cdots \bar{x}^{+}_{\alpha_{i}}(zq^{2m})=0\quad\textrm{for }m\geq c.\label{dfintegrabilnost}
\end{align}
Motivated by  \cite{G}, we defined in \cite{Kozic} (type 2) quantum quasi-particles of charge $m\in\mathbb{Z}_{>0}$ and color $i=1,2,...,n$ as coefficients of
the  operators
\begin{equation*}
\bar{x}_{m\alpha_{i}}^{+}(z)=\bar{x}^{+}_{\alpha_{i}}(z)\bar{x}^{+}_{\alpha_{i}}(zq^{2})\cdots\bar{x}_{\alpha_{i}}^{+}(zq^{2(m-1)})=\sum_{r\in\mathbb{Z}}\bar{x}^{+}_{m\alpha_{i}}(r)z^{-r-m}
\in\mathcal{E}(L).\end{equation*}

Fix $i=1,2,...,n$. Let $\mathcal{S}_i$ be a  set
$$\mathcal{S}_i =\left\{ \bar{x}_{\alpha_{i}}^{+}(z\q^l)\,\left|\, \right. l\in\mathbb{Z}_{\geq 0}\right\}\subset\mathcal{E}(L) .$$

\begin{pro}\label{partialthm}
There exists a unique smallest nonlocal $\q$-vertex
algebra $\left<\mathcal{S}_i \right>\subset \mathcal{E}(L)$ such that $\mathcal{S}_i\subset \left<\mathcal{S}_i\right>$.
\end{pro} 

\begin{prf}
Relation \eqref{dfkomutativnost} implies
$$\bar{x}_{ \alpha_{i}}^{+}(z_1 \q^{l_1})\bar{x}_{ \alpha_{i}}^{+}(z_2 \q^{l_2})\cdots
\bar{x}_{ \alpha_{i}}^{+}(z_m \q^{l_m}) \in\om(L,L((z_1,z_2,\ldots,z_m)))$$
for every positive integer $m$ and $l_1,...,l_m\geq 0$, so the statement
 follows directly from Theorem \ref{main}.
\end{prf}

\begin{rem}
For $i=1,2,...,n-1$ the pair $(\bar{x}_{\alpha_{i}}^{+}(z),\bar{x}_{\alpha_{i + 1}}^{+}(z))$  is not quasi-commutative, 
so  we can not use the results of  \Cref{teorijskisection} to construct a nonlocal $\q$-vertex
algebra containing $\bar{x}_{\alpha_{i}}^{+}(z)$ and $\bar{x}_{\alpha_{i+ 1}}^{+}(z)$.
\end{rem}

In the next two propositions we show that  nonlocal $\q$-vertex
algebra $\left<\mathcal{S}_i\right>$ provides an appropriate setting for studying
relations among quantum quasi-particles. Recall \eqref{intro3}. The following proposition is a direct consequence of \eqref{rtiprodukti} and \eqref{dfkomutativnost}.

\begin{pro}
For every positive integer $m$  we have
\begin{equation}\label{-1}
\bar{x}_{m\alpha_{i}}^{+}(z)=\bar{x}_{\alpha_{i}}^{+}(z)_{-1}\bar{x}_{\alpha_{i}}^{+}(z\q)_{-1}\ldots\bar{x}_{\alpha_{i}}^{+}(z\q^{m-1})_{-1}1\in \left<\mathcal{S}_i \right> .
\end{equation}
\end{pro}

Let $m$ and $k$ be positive integers such that
$m\leq k\leq c$.
First, we list $2m$ relations found in \cite{Kozic}:
\begin{align}
\bar{x}_{m\alpha_i}^{+}(z)\bar{x}_{k\alpha_i}^{+}(zq^{2m})=&\,\bar{x}_{(m+k)\alpha_i}^{+}(z),\tag{1}\label{11}\\
\bar{x}_{m\alpha_i}^{+}(zq^{2})\bar{x}_{k\alpha_i}^{+}(zq^{2m})=&\,\bar{x}_{\alpha_i}^{+}(zq^{2m})\bar{x}_{(m+k-1)\alpha_i}^{+}(zq^{2}),\tag{2}\label{2}\\
&\hspace{-6pt}\vdots\tag*{$\vdots\hspace{6pt}$}\\
\bar{x}_{m\alpha_i}^{+}(zq^{2(m-1)})\bar{x}_{k\alpha_i}^{+}(zq^{2m})=&\,\bar{x}_{(m-1)\alpha_i}^{+}(zq^{2m})\bar{x}_{(k+1)\alpha_i}^{+}(zq^{2(m-1)}),\tag{$m$}\label{mm}\\
\bar{x}_{m\alpha_i}^{+}(zq^{2k})\bar{x}_{k\alpha_i}^{+}(z)=&\,\bar{x}_{(m+k)\alpha_i}^{+}(z),\tag{$m+1$}\label{m+1}\\
\bar{x}_{m\alpha_i}^{+}(zq^{2(k-1)})\bar{x}_{k\alpha_i}^{+}(z)=&\,\bar{x}_{\alpha_i}^{+}(zq^{2(k-1)})\bar{x}_{(m+k-1)\alpha_i}^{+}(z),\tag{$m+2$}\label{m+2}\\
&\hspace{-6pt}\vdots\tag*{$\vdots\hspace{6pt}$}\\
\bar{x}_{m\alpha_i}^{+}(zq^{2(k-(m-1))})\bar{x}_{k\alpha_i}^{+}(z)=&\,\bar{x}_{(m-1)\alpha_i}^{+}(zq^{2(k-(m-1))})\bar{x}_{(k+1)\alpha_i}^{+}(z).\tag{$2m$}\label{2m}
\end{align}
It was proved in \cite{Kozic} that these relations are independent. More precisely, for any two integers $r,N$ and
a vector $v\in L$ we can express $2m$ vectors
\begin{equation}\label{vec1}
\bar{x}_{m\alpha_i}^{+}(r+l)\bar{x}_{k\alpha_i}^{+}(N-(r+l))v,\qquad
l=0,1,\ldots,2m-1,
\end{equation}
as a linear combination of  vectors
\begin{equation}\label{vec2}
\bar{x}_{m\alpha_i}^{+}(s)\bar{x}_{k\alpha_i}^{+}(N-s)v,\qquad
s\in\mathbb{Z}\setminus\left\{r,r+1,\ldots,r+2m-1\right\},
\end{equation}
 and vectors
\begin{equation}\label{vec3}
\bar{x}_{m'\alpha_i}^{+}(t)\bar{x}_{k'\alpha_i}^{+}(N-t)v,\qquad
0\leq m'<m,\quad m+k=m'+k',\quad t\in\mathbb{Z}.
\end{equation}
We shall say that a family $\mathcal{F}$  of equalities in
$\left<\mathcal{S}_i\right>$ is equivalent to equalities \eqref{11}--\eqref{2m} if for every $v\in L$ and
$r,N\in\mathbb{Z}$ we can, using only $\mathcal{F}$, express vectors
\eqref{vec1} as a linear combination of vectors \eqref{vec2} and \eqref{vec3}.

We  now formulate a vertex-operator theoretic interpretation of
\eqref{11}--\eqref{2m}, which is  analogous to relations \eqref{intro2} and \eqref{intro4}.
The proposition can be verified by a direct calculation.

\begin{pro}
For any two positive integers $m$ and $k$, $m\leq
k\leq c$, there exist scalars $c_{l,j},d_{l,j}\in\mathbb{C}(q)$, where
$j=0,1,\ldots,l$,  $l=0,1,\ldots,m-1$, such that  the following equalities hold for every $l=0,1,...,m-1$:
\begin{align}
&\bar{x}^{+}_{m\alpha_{i}}(z)_{-l-1}\bar{x}^{+}_{k\alpha_{i}}(z\textstyle\q)=
\displaystyle z^{-l}\sum_{j=0}^{l}c_{l,j}\bar{x}^{+}_{(m-1-j)\alpha_{i}}(z\textstyle\q^{l+1})\bar{x}^{+}_{(k+1+j)\alpha_{i}}(z\textstyle\q^{l-j})\in\left<\mathcal{S}_i\right>, \label{eq:1}\\
&\bar{x}^{+}_{m\alpha_{i}}(z\textstyle\q^k)_{-l-1}\bar{x}^{+}_{k\alpha_{i}}(z)=
\displaystyle z^{-l}\sum_{j=0}^{l}d_{l,j}\bar{x}^{+}_{j\alpha_{i}}(z\textstyle\q^{k+l-j})\bar{x}^{+}_{(m+k-j)\alpha_{i}}(z\textstyle\q^l)\in\left<\mathcal{S}_i\right>,\label{eq:2}
\end{align}
where $\bar{x}^{+}_{0\alpha_{i}}(z)=1$.
Furthermore, these $2m$ equalities in
$\left<\mathcal{S}_i\right>$ are equivalent to \eqref{11}--\eqref{2m}.
\end{pro}

\section*{Acknowledgement}
The author would like to acknowledge the support of the Australian Research Council
and of the Croatian Science Foundation under the project 2634.

The author would like to thank Mirko Primc for his valuable comments and suggestions.


\begin{thebibliography}{9}
\bibitem[A]{A1}
G. E. Andrews,
{\em The theory of partitions}, 
Encyclopedia of Mathematics and Its Applications, Vol. 2, Addison-Wesley, 1976.


\bibitem[AB]{AB}
I. I. Anguelova, M. J. Bergvelt,  
{\em $H_{D}$-Quantum vertex algebras and bicharacters},
Commun. Contemp. Math. \textbf{11} (2009) 937--991;
\href{http://arxiv.org/abs/0706.1528}{arXiv:0706.1528 [math.QA]}.

\bibitem[AKS]{AKS}
E. Ardonne, R. Kedem and M. Stone, 
{\em Fermionic characters of arbitrary highest-weight integrable $\widehat{sl}_{r+1}$-modules}, 
Comm. Math. Phys. \textbf{264} (2006), 427--464; 
\href{http://arxiv.org/abs/math/0504364}{arXiv:math/0504364 [math.RT]}.

\bibitem[BK]{Bak}
B. Bakalov, V. G. Kac,
{\em Field algebras}, 
Int. Math. Res. Not. (2003), no. 3, 123--159;
\href{http://arxiv.org/abs/math/0204282}{arXiv:math/0204282 [math.QA]}.

\bibitem[Ba]{Bar} 
I. Baranovi\'{c}, 
{\em Combinatorial bases of Feigin-Stoyanovsky's type subspaces of level 2 standard modules for $D_4^{(1)}$}, 
Comm. Algebra {\bf 39} (2011), 1007--1051;
\href{http://arxiv.org/abs/0903.0739}{arXiv:0903.0739 [math.QA]}.

\bibitem[B1]{Bor}
R. E. Borcherds, 
{\em Vertex algebras, Kac-Moody algebras, and the Monster}, 
Proc. Natl. Acad. Sci. USA \textbf{83} (1986), 3068--3071.

\bibitem[B2]{Bor2}
R. E. Borcherds, 
{\em Quantum vertex algebras},
in: 
Taniguchi Conference on Mathematics Nara'98, 
in:
Adv. Stud. Pure Math., 31, Math. Soc. Japan, Tokyo, 2001, 51--74.

\bibitem[Bu1]{Bu1} 
M. Butorac, 
{\em Combinatorial bases of principal subspaces for the affine Lie algebra of type $B_2^{(1)}$}, 
J. Pure Appl. Algebra \textbf{218} (2014), 424--447; 
\href{http://arxiv.org/abs/1212.5920}{arXiv:1212.5920 [math.QA]}.

\bibitem[Bu2]{Bu2} 
M. Butorac, 
{\em  Quasi-particle bases of principal subspaces for the affine Lie algebras of types $B_{l}^{(1)}$ and $C_{l}^{(1)}$},
Glas. Mat. Ser. III (to appear); 
\href{http://arxiv.org/abs/1505.00450}{arXiv:1505.00450 [math.QA]}.

\bibitem[CalLM1]{x1} 
C. Calinescu, J. Lepowsky and A. Milas, 
{\em Vertex-algebraic structure of the principal subspaces of certain $A^{(1)}_{1}$-modules, I: level one case}, 
Int. J. Math. \textbf{19}, no. 01, (2008), 71--92; 
\href{http://arxiv.org/abs/0704.1759}{arXiv:0704.1759 [math.QA]}.

\bibitem[CalLM2]{x2} 
C. Calinescu, J. Lepowsky and A. Milas, 
{\em Vertex-algebraic structure of the principal subspaces of certain $A^{(1)}_{1}$-modules, II: higher-level case}, 
J. Pure Appl. Algebra, \textbf{212} (2008), 1928--1950; 
\href{http://arxiv.org/abs/0710.1527}{arXiv:0710.1527 [math.QA]}.

\bibitem[CalLM3]{x3}
C. Calinescu, J. Lepowsky and A. Milas, 
{\em Vertex-algebraic structure of the principal subspaces of level one modules for the untwisted affine Lie algebras of types A,D,E}, 
J. Algebra \textbf{323} (2010), 167--192; 
\href{http://arxiv.org/abs/0908.4054}{arXiv:0908.4054 [math.QA]}.

\bibitem[CalLM4]{x4}
C. Calinescu, J. Lepowsky and A. Milas, 
{\em Vertex-algebraic structure of the principal subspaces ofstandard $A_2^{(2)}$-modules, I}, 
Int. J. Math. \textbf{25}, no. 07, 1450063 (2014);
\href{http://arxiv.org/abs/1402.3026}{arXiv:1402.3026 [math.QA]}.

\bibitem[CLM1]{CLM1}
S. Capparelli, J. Lepowsky and A. Milas,
{\em The Ro\-gers-Ra\-ma\-nu\-jan recur\-si\-on and intertwining o\-perators}, 
Comm. Con\-tem\-p. Math. \textbf{5} (2003), 947--966;
\href{http://arxiv.org/abs/math/0211265v2}{arXiv:math/0211265 [math.QA]}.

\bibitem[CLM2]{CLM2} 
S. Capparelli, J. Lepowsky and A. Milas, 
{\em The Rogers-Selberg recursions, the Gordon-Andrews identities and intertwining operators}, 
Ramanujan J. \textbf{12} (2006), 379--397; 
\href{http://arxiv.org/abs/math/0310080v2}{arXiv:math/0310080 [math.QA]}.


\bibitem[DF]{DF}
J. Ding, B. Feigin, 
{\em Commutative quantum current operators, semi-infinite construction and
functional models}, 
Represent. Theory 4, (2000), 330--341;
\href{http://arxiv.org/abs/q-alg/9612009}{arXiv:q-alg/9612009}.

\bibitem[DM]{DM}
J. Ding, T. Miwa, 
{\em Quantum current operators - I. Zeros and poles of quantum current operators and the condition of quantum integrability}, 
Publ. RIMS, Kyoto Univ. \textbf{33} (1997), 277--284;
\href{http://arxiv.org/abs/q-alg/9608001}{arXiv:q-alg/9608001}.

\bibitem[D]{D}
V. G. Drinfeld, 
{\em New realization of Yangian and quantized affine algebras}, 
Soviet Math. Dokl. 36 (1988), 212--216.

\bibitem[EK]{EK}
P. Etingof, D. Kazhdan,
{\em Quantization of Lie bialgebras V}, Selecta Math. (New Series) \textbf{6} (2000), 105--130;
\href{http://arxiv.org/abs/math/9808121}{arXiv:math/9808121 [math.QA]}.


\bibitem[FFJMM]{FFJMM} B. Feigin, E. Feigin, M. Jimbo, T. Miwa, E. Mukhin, 
{\em Principal $\widehat{sl_3}$ subspaces and quantum Toda Hamiltonian}, 
Adv. Stud. Pure Math. \textbf{54} (2009), 109--166;
\href{http://arxiv.org/abs/0707.1635}{arXiv:0707.1635 [math.QA]}.

\bibitem[FS]{FS}
A. V. Stoyanovsky, B. L. Feigin, 
{\em Functional models of the representations of current algebras, and semi-infinite Schubert cells}, (Russian) 
Funktsional. Anal. i Prilozhen. \textbf{28} (1994), no. 1, 68--90, 96; 
translation in Funct. Anal. Appl. \textbf{28} (1994), no. 1, 55--72; 
preprint B. L. Feigin and A. V. Stoyanovsky, {\em Quasi-particles models for the representations of Lie algebras and geometry of flag
manifold}, \href{http://arxiv.org/abs/hep-th/9308079}{arXiv:hep-th/9308079}.

\bibitem[FJ]{FJ}
I. B. Frenkel, N. Jing, 
{\em Vertex representations of quantum affine algebras}, 
Proc. Natl. Acad. Sci. USA, Vol. 85 (1988), 9373--9377.

\bibitem[FR]{FR}
E. Frenkel, N. Reshetikhin, 
{\em Towards deformed chiral algebras}, 
Quantum Group Symposium, XXI International Colloquium on Group Theoretical Methods in Physics (Goslar,
1996), Heron Press, Sofia, 1997, pp. 27--42;
\href{http://arxiv.org/abs/q-alg/9706023}{arXiv:q-alg/9706023}.

\bibitem[G]{G} 
G. Georgiev, 
{\em Combinatorial constructions of modules for infinite-dimensional Lie algebras, I. Principal subspace}, 
J. Pure Appl. Algebra \textbf{112} (1996), 247--286; 
\href{http://arxiv.org/abs/hep-th/9412054}{arXiv:hep-th/9412054}.

\bibitem[J1]{J1} 
M. Jerkovi\'{c}, 
{\em Recurrence relations for characters of affine Lie algebra $A_l^{(1)}$},  
J. Pure Appl. Algebra \textbf{213} (2009), 913--926;
\href{http://arxiv.org/abs/0803.1502}{arXiv:0803.1502 [math.QA]}. 

\bibitem[J2]{J2} 
M. Jerkovi\'{c}, 
{\em Character formulas for Feigin-Stoyanovsky's type subspaces of standard $\widetilde{\mathfrak{sl}}(3,C)$-modules}, 
Ramanujan J. \textbf{27} (2012), 357--376;
\href{http://arxiv.org/abs/1105.2927}{arXiv:1105.2927 [math.QA]}. 

\bibitem[JP]{JP}
M. Jerkovi\'{c}, M. Primc,
{\em Quasi-particle fermionic formulas for $(k,3)$-admissible configurations},
Cent. Eur. J. Math. \textbf{10} (2012), 703--721;
\href{http://arxiv.org/abs/1107.3900}{arXiv:1107.3900 [math.QA]}. 

\bibitem[K]{Kac} 
V. G. Kac, 
{\em Infinite-dimensional Lie algebras}, 
3rd ed., Cambridge University Press, Cambridge, 1990.

\bibitem[KC]{qKac}
V. Kac, P. Cheung,
{\em Quantum calculus},
Springer-Verlag, 2002.

\bibitem[Ka]{Kawa}
K. Kawasetsu,
{\em The Free Generalized Vertex Algebras and Generalized Principal Subspaces},
J. Algebra \textbf{444} (2015), 20--51;
\href{http://arxiv.org/abs/1502.05276}{arXiv:1502.05276 [math.QA]}.

\bibitem[Koy]{Koyama}
Y. Koyama, 
{\em Staggered Polarization of vertex models with $U_{q}(\tilde{\mathfrak{sl}}_{n})$-symmetry}, 
Comm. Math. Phys. \textbf{164}, no. 2 (1994), 277--291;
\href{http://arxiv.org/abs/hep-th/9307197}{arXiv:hep-th/9307197}.

\bibitem[Ko1]{Kozic}
S. Ko\v{z}i\'{c}, 
{\em Principal subspaces for quantum affine algebra $U_{q}(A_{n}^{(1)})$}, 
J. Pure Appl. Algebra \textbf{218} (2014), 2119--2148;
\href{http://arxiv.org/abs/1306.3712}{arXiv:1306.3712 [math.QA]}.

\bibitem[Ko2]{Kozic2}
S. Ko\v{z}i\'{c}, 
{\em A note on the zeroth products of Frenkel-Jing operators},
\href{http://arxiv.org/abs/1506.00050}{arXiv:1506.00050 [math.QA]}.

\bibitem[LP]{LP}
J. Lepowsky, M. Primc, 
{\em Structure of the standard modules for the affine Lie algebra
$A_{1}^{(1)}$}, Contemp. Math. 46, Amer. Math. Soc., Providence, 1985.

\bibitem[L1]{LiNonlocal}
H.-S. Li, 
{\em Nonlocal vertex algebras generated by formal vertex operators}, 
Selecta Math. (New Series) \textbf{11} (2005), 349--397;
\href{http://arxiv.org/abs/math/0502244}{arXiv:math/0502244 [math.QA]}.

\bibitem[L2]{Li2}
H.-S. Li, 
{\em Quantum vertex $\mathbb{F}((t))$-algebras and their modules}, 
J. Algebra \textbf{324} (2010), 2262--2304;
\href{http://arxiv.org/abs/0903.0186}{ 	arXiv:0903.0186 [math.QA]}.

\bibitem[L3]{Li3}
H.-S. Li, 
{\em Quantum vertex algebras and their $\phi$-coordinated modules},
Commun. Math. Phys. \textbf{308} (2011), 703--741;
\href{http://arxiv.org/pdf/0906.2710v2.pdf}{arXiv:0906.2710 [math.QA]}.

\bibitem[LL]{LiLep}
J. Lepowsky, H.-S. Li,
{\em Introduction to Vertex Operator Algebras and Their Representations},
Progress in Math., Vol. \textbf{227}, Birkhauser, Boston, 2004.

\bibitem[MP]{MiP} 
A. Milas, M. Penn, 
{\em Lattice vertex algebras and combinatorial bases: general case and W-algebras}, 
New York J. Math. \textbf{18} (2012), 621--650.

\bibitem[P1]{P1} 
M. Primc, 
{\em Vertex operator construction of standard modules for $A_n^{(1)}$}, 
Pacific J. Math \textbf{162} (1994), 143--187.

\bibitem[P2]{P2} 
M. Primc, 
{\em Basic Representations for classical affine Lie algebras}, 
J. Algebra \textbf{228} (2000), 1--50.


\bibitem[S1]{S1} 
C. Sadowski, 
{\em Presentations of the principal subspaces of the higher-level standard $\widehat{\mathfrak{sl}(3)}$-modules}, 
J. Pure Appl. Algebra \textbf{219} (2015) 2300-2345; 
\href{http://arxiv.org/abs/1312.6412}{arXiv:1312.6412 [math.QA]}.

\bibitem[S2]{S2} 
C. Sadowski, 
{\em Principal subspaces of higher-level standard $\widehat{\mathfrak{sl}(n)}$-modules},
Int. J. Math. \textbf{26}, no. 08, 1550053 (2015); 
\href{http://arxiv.org/abs/1406.0095}{arXiv:1406.0095 [math.QA]}.

\bibitem[T1]{T1} 
G. Trup\v{c}evi\'{c}, 
{\em Combinatorial bases of Feigin-Stoyanovsky's type subspaces of level 1 standard modules for $\widetilde{\mathfrak{sl}}(l+1,\mathbb{C})$}, 
Comm. Algebra \textbf{38} (2010), 3913--3940;
\href{http://arxiv.org/abs/0807.3363}{arXiv:0807.3363 [math.QA]}.

\bibitem[T2]{T2} 
G. Trup\v{c}evi\'{c}, 
{\em Combinatorial bases of Feigin-Stoyanovsky's type subspaces of higher-level standard $\widetilde{\mathfrak{sl}}(l+1,\mathbb{C})$-modules}, 
J. Algebra \textbf{322} (2009), 3744--3774;
\href{http://arxiv.org/abs/0810.5152}{arXiv:0810.5152 [math.QA]}.

\bibitem[T3]{T3}
G. Trup\v{c}evi\'{c}, 
{\em Characters of Feigin-Stoyanovsky's type subspaces of level one modules for affine Lie algebras of types $A_{l}^{(1)}$ and $D_4^{(1)}$}, 
Glas. Mat. Ser. III 46, \textbf{66} (2011), 49--70;
\href{http://arxiv.org/abs/1002.0348}{arXiv:1002.0348 [math.QA]}.

\end{thebibliography}
\end{document}